\documentclass[11pt]{article}
\usepackage{authblk}
\usepackage{amsmath,amssymb,amsthm,amsfonts,stmaryrd} 
\usepackage{tikz}
\usepackage[titletoc,toc,title]{appendix}
\usetikzlibrary{calc}
\usetikzlibrary{decorations.markings,arrows}
\usepackage{subfigure}
\usepackage{multirow}
\usepackage{epsf}
\usepackage[english]{babel}
\usepackage[latin1]{inputenc}
\usepackage[T1]{fontenc}
\usepackage{psfrag}
\usepackage{graphicx}
\newtheorem{theo}{Theorem}

\newtheorem{lemma}[theo]{Lemma}

\newtheorem{coro}[theo]{Corollary}
\newtheorem{obs}[theo]{Remark}

\begin{document}

\title{On a class of intersection graphs}

\author[1]{Mourad Ba\"{\i}ou}
\author[1] {Laurent Beaudou}
\author[2]{ Zhentao Li}
\author[1]{ Vincent Limouzy}

\affil[1]{CNRS and Universit\'e Clermont II, 
campus des c\'ezeaux \newline BP 125, 63173 Aubi\`ere cedex, France}

\affil[2]{CNRS and ENS Lyon}

\maketitle


\begin{abstract}
  Given a directed graph $D=(V,A)$ we define its intersection graph
  $I(D)=(A,E)$ to be the graph having $A$ as a node-set and two nodes
  of $I(D)$ are adjacent if their corresponding arcs share a common
  node that is the tail of at least one of these arcs. We call these
  graphs facility location graphs since they arise from the classical
  uncapacitated facility location problem. In this paper we
  show that facility location graphs are hard to recognize and they
  are easy to recognize when the underlying graph is triangle-free.
  We also determine the complexity of the vertex coloring, the stable
  set and the facility location problems on that class.
\end{abstract}




\section{Introduction}

\begin{figure}[htb]
  \begin{center}
    \begin{tikzpicture}[line cap=round,line join=round,x=.7cm,y=.7cm]
      \coordinate (a) at (0,0);
      \coordinate (b) at ($ (a) + (2,0) $);
      \coordinate (d) at ($ (a) + (6,0) $);
      \coordinate (e) at ($ (d) + (0,-1) $);
      \coordinate (c) at ($ (d) + (0,1) $);
      \coordinate (g) at ($ (d) + (2,0) $);
      \coordinate (h) at ($ (g) + (0,-1) $);
      \coordinate (f) at ($ (g) + (0,1) $);
      \coordinate (i) at ($ (g) + (3,.5) $);
      \coordinate (j) at ($ (i) + (-.75,-1) $);
      \coordinate (k) at ($ (i) + (.75,-1) $);
      \coordinate (l) at ($ (i) + (2,.5) $);
      \coordinate (m) at ($ (l) + (0,-2) $);
      \coordinate (n) at ($ (c) + (1,1) $);
      \coordinate (o) at ($ (i) + (0,1.5) $);
      \coordinate (p) at ($ (l) + (0,1) $);
      \coordinate (q) at ($ (a) + (1,-3) $);
      \coordinate (r) at ($ (i) + (0,-3.5) $);
      \coordinate (a1) at ($ (a) + (2,0.5) $);
      \coordinate (b1) at ($ (a1) + (-.75,-1) $);
      \coordinate (c1) at ($ (a1) + (+.75,-1) $);
      \coordinate (l1) at ($ (g) + (6,0) $);
      \coordinate (m1) at ($ (l1) + (-2,0) $);

      \fill (a1) circle (1.5pt);
      \fill (b1) circle (1.5pt);
      \fill (c1) circle (1.5pt);
      
      \fill (c) circle (1.5pt);
      \fill (d) circle (1.5pt);
      \fill (e) circle (1.5pt);
      \fill (f) circle (1.5pt);
      \fill (g) circle (1.5pt);
      \fill (h) circle (1.5pt);
      \fill (l) circle (1.5pt);
      \fill (m) circle (1.5pt);  
      
      \begin{scope}[decoration={
          markings,
          mark=at position 0.5 with {\arrow{stealth'}}}
        ] 
        
        \draw[postaction={decorate}] (e) -- (d) node[midway, left] {$a$};
        \draw[postaction={decorate}] (d) -- (c) node[midway, left] {$b$};
        
        \draw[postaction={decorate}] (f) -- (g) node[midway, left] {$b$};
        \draw[postaction={decorate}] (g) -- (h) node[midway, left] {$a$};
        
        \draw[postaction={decorate}] (a1) -- (b1)  node[midway, left] {$a$};
        \draw[postaction={decorate}] (a1) -- (c1)  node[midway,  right] {$b$};
        
        \draw[postaction={decorate}] (l) to[bend right] (m);  
        \draw (l1) node[left] {$a$}; 
        \draw (m1) node[right]  {$b$};  
        \draw[postaction={decorate}] (m) to[bend right] (l) ;
            
      \end{scope}    
 
    \end{tikzpicture}
  \end{center}
  \caption{The adjacency of two nodes $a$ and $b$ in $I(D)$.}
  \label{adjac}
\end{figure}
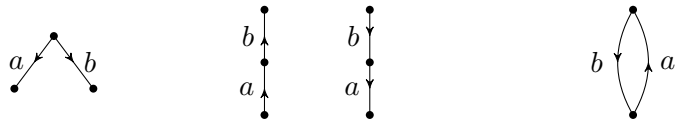

In this paper we study the following class of intersection
graphs. Given a directed graph $D=(V,A)$, we denote by $I(D)=(A,E)$
the {\it intersection graph of $D$} defined as follows:
\begin{itemize}
\item the node-set of $I(D)$ is the arc-set of $D$,
\item two nodes $a=(u,v)$ and $b=(w,t)$ of $I(D)$ are adjacent if one
  of the following holds: $u=w$ or $v=w$ or $t=u$ or $(u,v)=(t,w)$
  (see Figure \ref{adjac}).
\end{itemize}

We focus on two aspects: the recognition of these
intersection graphs and some combinatorial
optimization problem in this class. De Simone and Mannino
\cite{Simone} considered the recognition problem and provided a characterization of these graphs based on the structure of the (directed) neighborhood of a vertex. Unfortunately this characterization does not yield a polynomial time recognition algorithm.

Intersection graphs we consider arise from the {\it uncapacitated facility location
  problem} (UFLP) defined as follows. We are given a directed graph
$D=(V,A)$, costs $f(v)$ of opening a facility at node $v$ and cost $c(u,v)$ of assigning $v$ to $u$ (for each $(u,v) \in A$). We wish to select a subset of facilities to open and an assignment of each remaining nodes to a selected facility so as to minimize the cost of opening the selected facilities plus the cost of arcs used for assignment.

This problem can be formulated as a linear integer program as follows.

$$\mbox{min } \sum_{(u,v)\in A}c(u,v)x(u,v)+\sum_{v\in V}f(v)y(v)$$

$$\left\{ \begin{aligned}
\sum_{(u,v)\in A}x(u,v) + y(u) = 1 \quad & \forall u\in V,\\
x(u,v)\leq y(v) \quad & \forall (u,v)\in A,\\
x(u,v)\geq 0 \quad & \forall (u,v)\in A, \\
y(v)\geq 0 \quad & \forall v\in V,\\
x(u,v)\in \{0,1\} \quad & \forall (u,v)\in A,\\
y(v)\in \{0,1\} \quad & \forall v\in V.
\end{aligned} \right.$$

If we remove the variables $y(v)$ for all $v$ from the formulation
above, we get

$$\mbox{min } \sum_{(u,v)\in A}(c(u,v)-f(u))x(u,v)+\sum_{v\in V}f(v)$$

$$ \left\{ \begin{aligned}
\sum_{(u,v)\in A}x(u,v)  \leq 1 \quad  & \forall u\in V,\\
x(u,v)+\sum_{(v,w)\in A}x(v,w)  \leq 1 \quad & \forall (u,v)\in A,\\
x(u,v)\geq 0 \quad & \forall (u,v)\in A,\\
x(u,v)\in \{0,1\} \quad & \forall (u,v)\in A.
\end{aligned} \right.$$

This is exactly the maximal clique formulation of the {\it maximum
  stable set problem} associated with $I(D)$, where the weight of each
node $(u,v)$ of $I(D)$ is $f(u)-c(u,v)$. This correspondence is well
known in the literature (see in \cite{AVS, CorTh, Simone}). We may
consider several combinatorial optimization problems on directed graph
that may be reduce to the maximum stable set problem on an undirected
graph.  For example in \cite{Chv-Eben}, Chv\'atal and Ebenegger reduce
the max cut problem in a directed graph $D=(V,A)$ to the maximum
stable set problem in the following intersection graph called the {\it
  line graph of a directed graph}: we assign a node to each arc $a\in
A$ and two nodes are adjacent if the head of one (corresponding) arc is the tail of the other. They prove that recognizing such
graphs is \textsc{np}-complete. Balas \cite{Balas}
considered the asymmetric assignment problem. He defined an intersection graph
of a directed graph $D$ where nodes are arcs of $D$ and two nodes are adjacent if the two corresponding arcs have the same tail, the same head or the same extremities without being parallel. Balas uses this correspondence to develop new facets for the asymmetric assignment polytope.

We may generalize the notion of line graphs to directed graphs in many ways. The simplest involves deciding
\begin{enumerate}
\item
  if arcs that share a head are adjacent,
\item
  if arcs that share a tail are adjacent, and
\item
  if two arcs are adjacent if the head of one arc is the tail of the other. 
\end{enumerate}

It is not too difficult to show the recognition problem is easy if we choose non-adjacency for (3).

Choosing non-adjacency for (3) means that we could separate all vertices $v$ of a digraph $D$ into two vertices (one for all arcs entering that vertex and one for all arcs leaving it) and the line graph of the resulting digraph $D'$ is the same as the line graph of $D$. Furthermore, the line graph of $D'$ is the line graph of its underlying digraph. So all classes obtained by choosing non-adjacency for (3) are easy to recognize as it simply involves recognizing if a line graph is bipartite \cite{beineke_1970} (where some sides of the bipartition are possibly forced to have degree 1 from out choice of (1) and (2)).

So suppose arcs of type (3) are adjacent. Choosing adjacency for (1) and (2) gives the line graphs of the underlying undirected graph, and these are easy to recognize \cite{beineke_1970}. Choosing non-adjacency for both (1) and (2) leads to the line graphs defined by Chv\'atal and Ebenegger and it is \textsc{np}-complete to recognize them \cite{Chv-Eben}. And picking exactly one of (1) and (2) to be adjacent and non-adjacency for the other leads to the same class of graphs (as we can simply reverse all arcs of a digraph before taking its line graph) and we wish to determine the complexity of recognizing this very last class.


Finally, note that since the stable set problem in our class is equivalent to the facility location problem, one may use all the material
developed for facility location problem to solve the stable set
problem in these graphs.  It is well known that in practice the
facility location problem may be solved efficiently via several
approaches: polyhedra, approximation algorithms and heuristics.

This paper is organized as follows. Section \ref{def} contains some basic definition and notations. Other definitions and notations will
be given when needed. In section \ref{RLFG}, we show that
facility location graphs are hard to recognize and in
Section \ref{RTFLG} we show that the subclass of triangle-free
facility location graphs are recognizable in polynomial time. Section
\ref{RP} is devoted to some combinatorial optimization problem in
facility location graphs. In particlular we show that the maximum
stable set problem remains \textsc{np}-complete in triangle-free facility
location graphs but the vertex coloring problem is solvable in polynomial time in this class. We also discuss the facility location problem and show it is
\textsc{np}-complete in some restricted class of graphs. We provide
concluding remarks in Section \ref{conc}.

\section{Definitions and notations}

\label{def}

Let $G$ be an undirected graph, we say that $G$ is a {\it facility
  location} (FL) graph if there exists a directed graph $D$ such that
$G=I(D)$. Any FL graph will be denoted by $I(D)$, this notation helps
to indicate the directed graph $D$ from which our FL graph may be
obtained. Such a graph $D$ is called the {\em preimage} of $G$.

Let $D=(V,A)$ be a directed graph. Given an arc $a=(u,v)\in A$, the
node $u$ is called the {\it tail} of $a$ and $v$ is called the {\it
  head} of $a$. Sometimes we use the notation $t(a)$ (respectively
$h(a)$) for the tail (respectively the head) of $a$.  A {\it sink} is
a node which is a tail of no arc in $A$. A {\it branch} in a directed
graph is an arc $(u,v)$ where $v$ is a sink and is the head of only
$(u,v)$.

A {\it cycle} $C$ in $D$ is a cycle in the underlying undirected graph of $D$. I.e., an ordered sequence of arcs $a_1, a_2,
\dots,a_p,a_{p+1}$, where $a_i$ and $a_{i+1}$ are incident, for
$i=1,\dots,p$, with $a_{p+1}=a_1$. If $a_1$ and $a_p$ are not
incident, then this sequence is called a {\it path}.  We denote by $A(C)$ the arcs of
$C$ and by $V(C)$ its nodes, that are the endnodes of the arcs in
$A(C)$. The nodes of $V(C)$ may be partitioned into three sets (1)
$\dot{C}$, the nodes that are the tail of two arcs in $A(C)$, (2)
$\hat{C}$, the nodes that are the head of two arcs in $A(C)$ and (3)
$\tilde{C}$, the nodes that are the tail of one arc and the head of
the other arc in $A(C)$. When $\hat{C}$ (or $\dot{C}$) is empty, the
cycle $C$ is the classical {\it directed} cycle. Similarly the nodes
of a path $P$, except its extremities, may be partitioned into three
sets $\dot{P}$, $\hat{P}$ and $\tilde{P}$.  When
$\hat{P}=\dot{P}=\emptyset$, $P$ is called a {\it directed path}.  We
define a cycle $C$ in an undirected graph as a sequence of ordered
nodes, instead of ordered edges, this is useful when we study the
correspondence between $I(D)$ and $D$. The nodes of $C$ are denoted by
$A(C)$ and its edges by $E(C)$. A {\it path} is defined similarly.

Define $x_1,\dots,x_n$ to be $n$ Boolean variables. Define a {\em literal} $\lambda_i$ be either a Boolean variable $x_i$ or its negation $\bar{x_i}$. A {\it clause} $C$, is a disjunction of literals $\lambda_i$, that is $C=\lambda_{i_1}\vee\dots\vee\lambda_{i_k}$. $F=C_1\wedge\dots\wedge C_m$ is a conjunction of $m$ clauses. In the {\it satisfiablity problem} \textsc{sat}, we want to decide if there exists values $x_i$ such that an input conjunction of disjunctions (of literals) $F$ evaluates to {\em true}. If such values exist, we say $F$ is {\em satisfiable}. If each clause is a disjunction of at most three literals then the problem is
called 3-satisfiability. The problem 3-\textsc{sat} has been shown
\textsc{np}-complete by Karp \cite{karp}.

An undirected graph $G$ is triangle-free if it does not contain a
clique of size 3.  A {\it wheel} $W_n$ is a graph obtained from a cycle $C_n$by adding a vertex adjacent to all vertices of the cycle.

\section{Recognizing facility location graphs is \textsc{np}-complete}
\label{RLFG}

The main result of this section is the following: 

\begin{theo}
\label{NP-complete}
Recognizing facility location graphs is \textsc{np}-complete.
\end{theo}

The proof of this theorem is given in subsection
\ref{proof-NP-complete}.  We first give a sketch of this proof and
some useful lemmas before providing the detailed proof.

\subsection{Proof sketch}
We will reduce the problem 3-\textsc{sat} to the recognition of FL
graphs.  We assume we are given an instance of the problem
3-\textsc{sat}. That is, we have $n$ Boolean variables $x_1,\dots,x_n$
and a Boolean formula $F=C_1\wedge\dots\wedge C_m$, where each clause
$C_j=\lambda_{j_1}\vee\lambda_{j_2}\vee\lambda_{j_3}$, for
$j=1\dots,m$.  From $F$ we construct an undirected graph $G_F$ and we
show that $F$ is satisfiable if and only if $G_F$ is a facility location
graph.
 
We build $G_F$ using gadgets for variables and clauses. Values for variables are stored, replicated and negated through the ``branches'' of the variable gadgets. These branches are then connected to the clauses gadgets of clauses that contain these variables (and their negation). 

More precisely, the construction of $G_F$ follows three steps: (1) for each
variable $x_i$, we construct a graph called $\textsc{Gad}^1_i$
(\textsc{Gad} stands for gadget), (2) for each clause $C_j$, another
gadget called $\textsc{Gad}^2_j$ is constructed and (3) we connect the
graphs $\textsc{Gad}^1_i$ and $\textsc{Gad}^2_j$ to produce
$G_F$. Each graph $\textsc{Gad}^1_i$ contains $2m$ branches where each
branch express the fact that the variable $x_i$ (or $\bar{x}_i$) is
present in the clause $C_j$, $j=1,\dots,m$. Each graph
$\textsc{Gad}^2_j$ contains exactly three branches where each branch
expresses the litterals of this clause $\lambda_{j_1}$, $\lambda_{j_2}$
and $\lambda_{j_3}$.

The three following subsections are devoted to the construction of the
graphs $\textsc{Gad}^1_i$, $\textsc{Gad}^2_j$ and $G_F$.

\subsection{The construction of  the graphs $\textsc{Gad}^1_i$}
 
\begin{obs} 
\label{wheel-obs}
There are 15 directed graphs whose intersection graph is the wheel
$W_5$. We list 5 of them that will be useful for our reduction. From
these graphs, we obtain all the remaining directed graphs by
identifying the head of the pending arc with the tail of one of the
arcs entering this pending arc.
\end{obs}

 \begin{figure}[htb]
    \begin{center}
      \begin{tikzpicture}[line cap=round,line join=round,x=1cm,y=1cm]
        \coordinate (a) at (-2,1.5);
        \coordinate (b) at ($ (a) + (126:1.5) $);
        \coordinate (c) at ($ (a) + (-162:1.5) $);
        \coordinate (d) at ($ (a) + (-90:1.5) $);
        \coordinate (e) at ($ (a) + (-18:1.5) $);
        \coordinate (f) at ($ (a) + (54:1.5) $);
        
        \draw (b) -- (c) -- (d) -- (e) -- (f) -- (b);
        \draw (a) -- (b);
        \draw (a) -- (c);
        \draw (a) -- (d);
        \draw (a) -- (e);
        \draw (a) -- (f);

        \fill (a) circle (1.5pt);
        \fill (b) circle (1.5pt);
        \fill (c) circle (1.5pt);
        \fill (d) circle (1.5pt);
        \fill (e) circle (1.5pt);
        \fill (f) circle (1.5pt);
        
        \draw (a) node[above] {$a$};
        \draw (b) node[left] {$b$};
        \draw (c) node[left] {$c$};
        \draw (d) node[below] {$d$};
        \draw (e) node[right] {$e$};
        \draw (f) node[right] {$f$};

 \end{tikzpicture}
    \end{center}
    \caption{The wheel $W_5$.}
    \label{fig:w5}
  \end{figure}
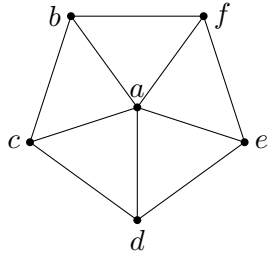
  
   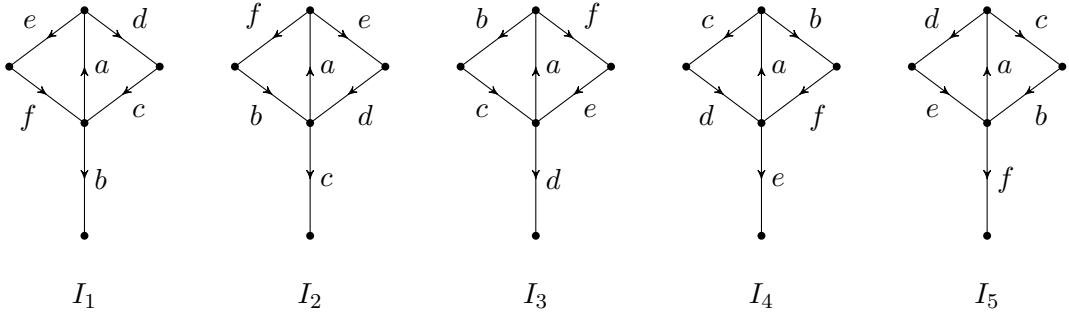
\begin{figure}[htb]
    \begin{center}
      \begin{tikzpicture}[line cap=round,line join=round,x=1cm,y=1cm]

   \coordinate (x) at (3,2.25);
        \coordinate (y) at (2,3);
        \coordinate (z) at (1,2.25);
        \coordinate (t) at (2,1.5);
        \coordinate (u) at (2,0);
       \coordinate (v) at (2,-0.5);
        
        \begin{scope}[decoration={
            markings,
            mark=at position 0.5 with {\arrow{stealth'}}}
          ] 
          
          \draw[postaction={decorate}] (y) -- (z) node[midway,above left]{$d$};
          \draw[postaction={decorate}] (y) -- (x) node[midway,above right]{$c$};
          \draw[postaction={decorate}] (z) -- (t) node[midway,below left]{$e$};
          \draw[postaction={decorate}] (x) -- (t) node[midway,below right]{$b$};
          \draw[postaction={decorate}] (t) -- (y) node[midway,right]{$a$};
          \draw[postaction={decorate}] (t) -- (u) node[midway,right]{$f$};   
           \draw (v) node[below] {$I_5$};
          
        \fill (x) circle (1.5pt);
        \fill (y) circle (1.5pt);
        \fill (z) circle (1.5pt);
        \fill (t) circle (1.5pt);
        \fill (u) circle (1.5pt);
                 
        \end{scope}

   \coordinate (x) at (0,2.25);
        \coordinate (y) at (-1,3);
        \coordinate (z) at (-2,2.25);
        \coordinate (t) at (-1,1.5);
        \coordinate (u) at (-1,0);
        \coordinate (v) at (-1,-0.5);
  
        \begin{scope}[decoration={
            markings,
            mark=at position 0.5 with {\arrow{stealth'}}}
          ] 
          
          \draw[postaction={decorate}] (y) -- (z) node[midway,above left]{$c$};
          \draw[postaction={decorate}] (y) -- (x) node[midway,above right]{$b$};
          \draw[postaction={decorate}] (z) -- (t) node[midway,below left]{$d$};
          \draw[postaction={decorate}] (x) -- (t) node[midway,below right]{$f$};
          \draw[postaction={decorate}] (t) -- (y) node[midway,right]{$a$};
          \draw[postaction={decorate}] (t) -- (u) node[midway,right]{$e$};   
          \draw (v) node[below] {$I_4$};
          
        \fill (x) circle (1.5pt);
        \fill (y) circle (1.5pt);
        \fill (z) circle (1.5pt);
        \fill (t) circle (1.5pt);
        \fill (u) circle (1.5pt);
                 
        \end{scope}

        \coordinate (x) at (-3,2.25);
        \coordinate (y) at (-4,3);
        \coordinate (z) at (-5,2.25);
        \coordinate (t) at (-4,1.5);
        \coordinate (u) at (-4,0);
        \coordinate (v) at (-4,-0.5);
        
        \begin{scope}[decoration={
            markings,
            mark=at position 0.5 with {\arrow{stealth'}}}
          ] 
          
          \draw[postaction={decorate}] (y) -- (z) node[midway,above left]{$b$};
          \draw[postaction={decorate}] (y) -- (x) node[midway,above right]{$f$};
          \draw[postaction={decorate}] (z) -- (t) node[midway,below left]{$c$};
          \draw[postaction={decorate}] (x) -- (t) node[midway,below right]{$e$};
          \draw[postaction={decorate}] (t) -- (y) node[midway,right]{$a$};
          \draw[postaction={decorate}] (t) -- (u) node[midway,right]{$d$};   
           \draw (v) node[below] {$I_3$};
          
        \fill (x) circle (1.5pt);
        \fill (y) circle (1.5pt);
        \fill (z) circle (1.5pt);
        \fill (t) circle (1.5pt);
        \fill (u) circle (1.5pt);
                 
        \end{scope}
 
         \coordinate (x) at (-6,2.25);
        \coordinate (y) at (-7,3);
        \coordinate (z) at (-8,2.25);
        \coordinate (t) at (-7,1.5);
        \coordinate (u) at (-7,0);
        \coordinate (v) at (-7,-0.5);
        
        \begin{scope}[decoration={
            markings,
            mark=at position 0.5 with {\arrow{stealth'}}}
          ] 
          
          \draw[postaction={decorate}] (y) -- (z) node[midway,above left]{$f$};
          \draw[postaction={decorate}] (y) -- (x) node[midway,above right]{$e$};
          \draw[postaction={decorate}] (z) -- (t) node[midway,below left]{$b$};
          \draw[postaction={decorate}] (x) -- (t) node[midway,below right]{$d$};
          \draw[postaction={decorate}] (t) -- (y) node[midway,right]{$a$};
          \draw[postaction={decorate}] (t) -- (u) node[midway,right]{$c$};   
           \draw (v) node[below] {$I_2$};
          
        \fill (x) circle (1.5pt);
        \fill (y) circle (1.5pt);
        \fill (z) circle (1.5pt);
        \fill (t) circle (1.5pt);
        \fill (u) circle (1.5pt);
                 
        \end{scope}

          \coordinate (x) at (-9,2.25);
        \coordinate (y) at (-10,3);
        \coordinate (z) at (-11,2.25);
        \coordinate (t) at (-10,1.5);
        \coordinate (u) at (-10,0);
        \coordinate (v) at (-10,-0.5);
        
        \begin{scope}[decoration={
            markings,
            mark=at position 0.5 with {\arrow{stealth'}}}
          ] 
          
          \draw[postaction={decorate}] (y) -- (z) node[midway,above left]{$e$};
          \draw[postaction={decorate}] (y) -- (x) node[midway,above right]{$d$};
          \draw[postaction={decorate}] (z) -- (t) node[midway,below left]{$f$};
          \draw[postaction={decorate}] (x) -- (t) node[midway,below right]{$c$};
          \draw[postaction={decorate}] (t) -- (y) node[midway,right]{$a$};
          \draw[postaction={decorate}] (t) -- (u) node[midway,right]{$b$};   
          \draw (v) node[below] {$I_1$};
          
        \fill (x) circle (1.5pt);
        \fill (y) circle (1.5pt);
        \fill (z) circle (1.5pt);
        \fill (t) circle (1.5pt);
        \fill (u) circle (1.5pt);
                 
        \end{scope}

      \end{tikzpicture}
    \end{center}
    \caption{The 5 oriented graphs having $W_5$ as their intersection graph.}
    \label{fig:w5isflg}
  \end{figure}

\begin{figure}[htb]
    \begin{center}
      \begin{tikzpicture}[line cap=round,line join=round,x=1cm,y=1cm]
        \coordinate (a) at (-2,1.5);
        \coordinate (b) at ($ (a) + (126:1.5) $);
        \coordinate (c) at ($ (a) + (-162:1.5) $);
        \coordinate (d) at ($ (a) + (-90:1.5) $);
        \coordinate (e) at ($ (a) + (-18:1.5) $);
        \coordinate (f) at ($ (a) + (54:1.5) $);
        \coordinate (g) at ($ (b) + (90:1.5) $);
        \coordinate (h) at ($ (f) + (90:1.5) $);
        \coordinate (i) at ($ (d) + (-2,0) $);
        \coordinate (j) at ($ (d) + (2,0) $);
         \coordinate (x) at ($(d) + (0,-1)$);
        
        \draw (b) -- (c) -- (d) -- (e) -- (f) -- (b);
        \draw (c) -- (i) -- (d) -- (c);
        \draw (d) -- (j) -- (e) -- (d);
        \draw (a) -- (b);
        \draw (a) -- (c);
        \draw (a) -- (d);
        \draw (a) -- (e);
        \draw (a) -- (f);
        \draw (b) -- (g);
        \draw (f) -- (h);
         \draw (x) node[below] {$I$};

        \fill (a) circle (1.5pt);
        \fill (b) circle (1.5pt);
        \fill (c) circle (1.5pt);
        \fill (d) circle (1.5pt);
        \fill (e) circle (1.5pt);
        \fill (f) circle (1.5pt);
        \fill (g) circle (1.5pt);
        \fill (h) circle (1.5pt);
        \fill (i) circle (1.5pt);
        \fill (j) circle (1.5pt);
       
        \draw (a) node[above] {$a$};
        \draw (b) node[left] {$b$};
        \draw (c) node[left] {$c$};
        \draw (d) node[below] {$d$};
        \draw (e) node[right] {$e$};
        \draw (f) node[right] {$f$}; 
        \draw (g) node[above] {$g$};
        \draw (h) node[above] {$h$}; 
        \draw (i) node[left] {$i$};
        \draw (j) node[right] {$j$};

       \coordinate (a) at (4,1.5);
        \coordinate (b) at ($ (a) + (126:1.5) $);
        \coordinate (c) at ($ (a) + (-162:1.5) $);
        \coordinate (d) at ($ (a) + (-90:1.5) $);
        \coordinate (e) at ($ (a) + (-18:1.5) $);
        \coordinate (f) at ($ (a) + (54:1.5) $);
        \coordinate (g) at ($ (b) + (90:1.5) $);
        \coordinate (h) at ($ (f) + (90:1.5) $);
        \coordinate (i) at ($ (d) + (-2,0) $);
        \coordinate (j) at ($ (d) + (2,0) $);
         \coordinate (x) at ($(d) + (0,-1)$);
        
        \draw (b) -- (c) -- (d) -- (e) -- (f) -- (b);
        \draw (c) -- (i) -- (d) -- (c);
        \draw (d) -- (j) -- (e) -- (d);
        \draw (a) -- (b);
        \draw (a) -- (c);
        \draw (a) -- (d);
        \draw (a) -- (e);
        \draw (a) -- (f);
        \draw (b) -- (g);
        \draw (f) -- (h);
         \draw (x) node[below] {$I'$};

        \fill (a) circle (1.5pt);
        \fill (b) circle (1.5pt);
        \fill (c) circle (1.5pt);
        \fill (d) circle (1.5pt);
        \fill (e) circle (1.5pt);
        \fill (f) circle (1.5pt);
        \fill (g) circle (1.5pt);
        \fill (h) circle (1.5pt);
        \fill (i) circle (1.5pt);
        \fill (j) circle (1.5pt);
       
        \draw (a) node[above] {$a'$};
        \draw (b) node[left] {$b'$};
        \draw (c) node[left] {$c'$};
        \draw (d) node[below] {$d'$};
        \draw (e) node[right] {$e'$};
        \draw (f) node[right] {$f'$}; 
        \draw (g) node[above] {$g'$};
        \draw (h) node[above] {$h'$}; 
        \draw (i) node[left] {$j'$};
        \draw (j) node[right] {$i'$};

              \end{tikzpicture}
    \end{center}
    \caption{$I$ and $I'$ are extensions of $W_5$.}
    \label{fig:demi_inv}
  \end{figure}
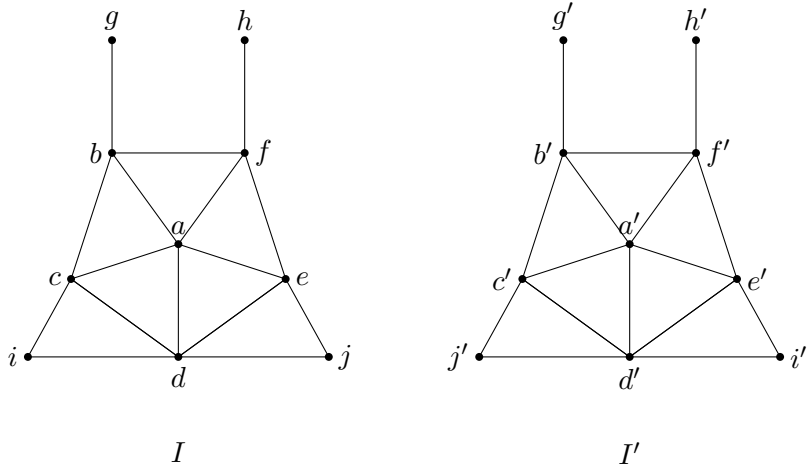

  Consider the two undirected graphs $I$ and $I'$ of Figure
  \ref{fig:demi_inv}. These two graphs are isomorphic, they are
  obtained from the wheel $W_5$ by adding four nodes. This restricts
  the number of possible preimage for the wheel to only 2. When we
  join the nodes $j$ and $j'$ we obtain the desired graph called
  \textsc{Inv} which is again the intersection graph of exactly two
  directed graphs. The graph \textsc{Inv} will be useful for the
  construction of $\textsc{Gad}^2_j$.

  Recall that $W_5$ may be the intersection of 15 directed graphs as
  described by Remark \ref{wheel-obs}. Let us discuss the possible
  directed graphs having $I$ as their intersection graph. The subgraph
  of $I$ induced by the nodes $a,b,c,d,e,f$ is a wheel $W_5$ and thus
  is the intersection of 15 directed graphs obtained from those of
  Figure \ref{fig:w5isflg} as described by Remark \ref{wheel-obs}. In
  $I$ there are two pendent nodes $g$ and $h$ adjacent to two
  neighbors $b$ and $f$ in $I$. This implies that the adjacency
  between $b$ and $f$ cannot be of the form as represented in the
  graphs $I_2$, $I_3$ and $I_4$ of Figure \ref{fig:w5isflg}. Hence we
  have only two possible directed graphs $I_1$ or $I_5$, since in
  neither of these two graphs we may identify the head of the pending
  arc with the tail of one of the arcs entering this pending
  arc. Moreover, there is one way of adding the arcs $i$ and $j$ in
  each case. Consequently there are only two graphs $D_1$ and $D_2$
  such that $I(D_1)=I(D_2)=I$. These graphs are represented in Figure
  \ref{fig:demi_inv1}.  It is obvious that there are also two directed
  graphs whose intersection graph is $I'$. For convenience these two
  graphs will be shown in Figure \ref{fig:demi_inv2}.
  
 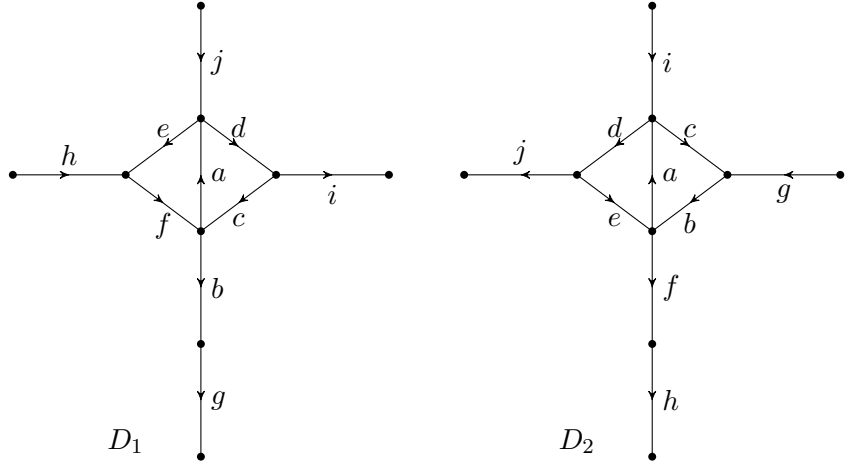
\begin{figure}[htb]
    \begin{center}
      \begin{tikzpicture}[line cap=round,line join=round,x=1cm,y=1cm]

        \coordinate (x) at (4,2.25);
        \coordinate (y) at (3,3);
        \coordinate (z) at (2,2.25);
        \coordinate (t) at (3,1.5);
        \coordinate (u) at (3,0);
        \coordinate (v) at (3,-1.5);
        \coordinate (p) at (3,4.5);
        \coordinate (q) at (0.5,2.25);
        \coordinate (w) at (5.5,2.25);
        \coordinate (k) at (2,-1);
        
        \begin{scope}[decoration={
            markings,
            mark=at position 0.5 with {\arrow{stealth'}}}
          ] 
        
        \draw[postaction={decorate}] (y) -- (z) node[midway,left, above]{$e$};
        \draw[postaction={decorate}] (y) -- (x) node[midway,right, above]{$d$};
        \draw[postaction={decorate}] (z) -- (t) node[midway,left, below]{$f$};
        \draw[postaction={decorate}] (x) -- (t) node[midway,right, below]{$c$};
        \draw[postaction={decorate}] (t) -- (y) node[midway,right]{$a$};
        \draw[postaction={decorate}] (t) -- (u) node[midway,right]{$b$};
        \draw[postaction={decorate}] (p) -- (y) node[midway,right]{$j$};
        \draw[postaction={decorate}] (x) -- (w) node[midway, below]{$i$};
        \draw[postaction={decorate}] (q) -- (z) node[midway,above]{$h$};
        \draw[postaction={decorate}] (u) -- (v) node[midway,right]{$g$};
        \draw (k) node[below] {$D_1$};
      \end{scope}

         \fill (x) circle (1.5pt);
        \fill (y) circle (1.5pt);
        \fill (z) circle (1.5pt);
        \fill (t) circle (1.5pt);
        \fill (u) circle (1.5pt);
        \fill (p) circle (1.5pt);
        \fill (q) circle (1.5pt);
        \fill (v) circle (1.5pt);
        \fill (w) circle (1.5pt);

         \coordinate (x) at (10,2.25);
        \coordinate (y) at (9,3);
        \coordinate (z) at (8,2.25);
        \coordinate (t) at (9,1.5);
        \coordinate (u) at (9,0);
        \coordinate (v) at (9,-1.5);
        \coordinate (p) at (9,4.5);
        \coordinate (q) at (6.5,2.25);
        \coordinate (w) at (11.5,2.25);
         \coordinate (k) at (8,-1);
         
        \begin{scope}[decoration={
            markings,
            mark=at position 0.5 with {\arrow{stealth'}}}
          ] 
        
        \draw[postaction={decorate}] (y) -- (z) node[midway,left, above]{$d$};
        \draw[postaction={decorate}] (y) -- (x) node[midway,right, above]{$c$};
        \draw[postaction={decorate}] (z) -- (t) node[midway,left, below]{$e$};
        \draw[postaction={decorate}] (x) -- (t) node[midway,right, below]{$b$};
        \draw[postaction={decorate}] (t) -- (y) node[midway,right]{$a$};
        \draw[postaction={decorate}] (t) -- (u) node[midway,right]{$f$};
        \draw[postaction={decorate}] (p) -- (y) node[midway,right]{$i$};
        \draw[postaction={decorate}] (w) -- (x) node[midway, below]{$g$};
        \draw[postaction={decorate}] (z) -- (q) node[midway,above]{$j$};
        \draw[postaction={decorate}] (u) -- (v) node[midway,right]{$h$};
         \draw (k) node[below] {$D_2$};
      \end{scope}

         \fill (x) circle (1.5pt);
        \fill (y) circle (1.5pt);
        \fill (z) circle (1.5pt);
        \fill (t) circle (1.5pt);
        \fill (u) circle (1.5pt);
        \fill (p) circle (1.5pt);
        \fill (q) circle (1.5pt);
        \fill (v) circle (1.5pt);
        \fill (w) circle (1.5pt);

      \end{tikzpicture}
    \end{center}
    \caption{The two directed graphs having $I$ as an intersection graph.}
    \label{fig:demi_inv1}
  \end{figure}
   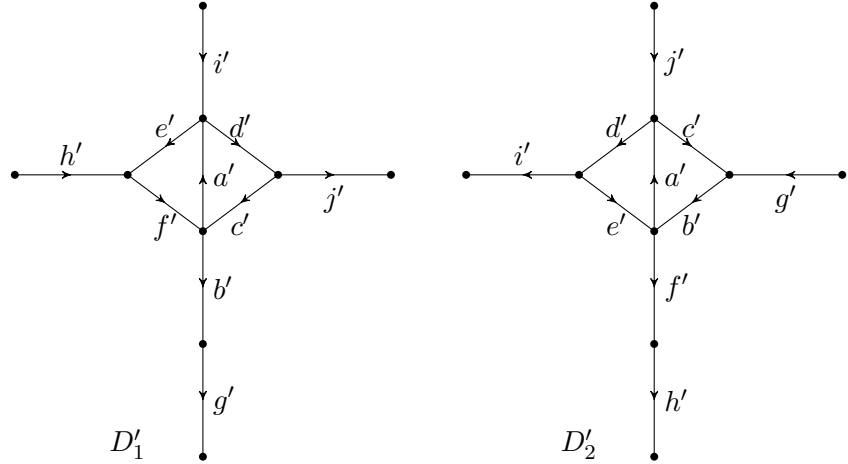
\begin{figure}[htb]
    \begin{center}
      \begin{tikzpicture}[line cap=round,line join=round,x=1cm,y=1cm]

        \coordinate (x) at (4,2.25);
        \coordinate (y) at (3,3);
        \coordinate (z) at (2,2.25);
        \coordinate (t) at (3,1.5);
        \coordinate (u) at (3,0);
        \coordinate (v) at (3,-1.5);
        \coordinate (p) at (3,4.5);
        \coordinate (q) at (0.5,2.25);
        \coordinate (w) at (5.5,2.25);
        \coordinate (k) at (2,-1);
        
        \begin{scope}[decoration={
            markings,
            mark=at position 0.5 with {\arrow{stealth'}}}
          ] 
        
        \draw[postaction={decorate}] (y) -- (z) node[midway,left, above]{$e'$};
        \draw[postaction={decorate}] (y) -- (x) node[midway,right, above]{$d'$};
        \draw[postaction={decorate}] (z) -- (t) node[midway,left, below]{$f'$};
        \draw[postaction={decorate}] (x) -- (t) node[midway,right, below]{$c'$};
        \draw[postaction={decorate}] (t) -- (y) node[midway,right]{$a'$};
        \draw[postaction={decorate}] (t) -- (u) node[midway,right]{$b'$};
        \draw[postaction={decorate}] (p) -- (y) node[midway,right]{$i'$};
        \draw[postaction={decorate}] (x) -- (w) node[midway, below]{$j'$};
        \draw[postaction={decorate}] (q) -- (z) node[midway,above]{$h'$};
        \draw[postaction={decorate}] (u) -- (v) node[midway,right]{$g'$};
        \draw (k) node[below] {$D'_1$};
      \end{scope}

         \fill (x) circle (1.5pt);
        \fill (y) circle (1.5pt);
        \fill (z) circle (1.5pt);
        \fill (t) circle (1.5pt);
        \fill (u) circle (1.5pt);
        \fill (p) circle (1.5pt);
        \fill (q) circle (1.5pt);
        \fill (v) circle (1.5pt);
        \fill (w) circle (1.5pt);

         \coordinate (x) at (10,2.25);
        \coordinate (y) at (9,3);
        \coordinate (z) at (8,2.25);
        \coordinate (t) at (9,1.5);
        \coordinate (u) at (9,0);
        \coordinate (v) at (9,-1.5);
        \coordinate (p) at (9,4.5);
        \coordinate (q) at (6.5,2.25);
        \coordinate (w) at (11.5,2.25);
         \coordinate (k) at (8,-1);
         
        \begin{scope}[decoration={
            markings,
            mark=at position 0.5 with {\arrow{stealth'}}}
          ] 
        
        \draw[postaction={decorate}] (y) -- (z) node[midway,left, above]{$d'$};
        \draw[postaction={decorate}] (y) -- (x) node[midway,right, above]{$c'$};
        \draw[postaction={decorate}] (z) -- (t) node[midway,left, below]{$e'$};
        \draw[postaction={decorate}] (x) -- (t) node[midway,right, below]{$b'$};
        \draw[postaction={decorate}] (t) -- (y) node[midway,right]{$a'$};
        \draw[postaction={decorate}] (t) -- (u) node[midway,right]{$f'$};
        \draw[postaction={decorate}] (p) -- (y) node[midway,right]{$j'$};
        \draw[postaction={decorate}] (w) -- (x) node[midway, below]{$g'$};
        \draw[postaction={decorate}] (z) -- (q) node[midway,above]{$i'$};
        \draw[postaction={decorate}] (u) -- (v) node[midway,right]{$h'$};
         \draw (k) node[below] {$D'_2$};
      \end{scope}

         \fill (x) circle (1.5pt);
        \fill (y) circle (1.5pt);
        \fill (z) circle (1.5pt);
        \fill (t) circle (1.5pt);
        \fill (u) circle (1.5pt);
        \fill (p) circle (1.5pt);
        \fill (q) circle (1.5pt);
        \fill (v) circle (1.5pt);
        \fill (w) circle (1.5pt);

      \end{tikzpicture}
    \end{center}
    \caption{The two directed graphs having $I'$ as an intersection graph.}
    \label{fig:demi_inv2}
  \end{figure}

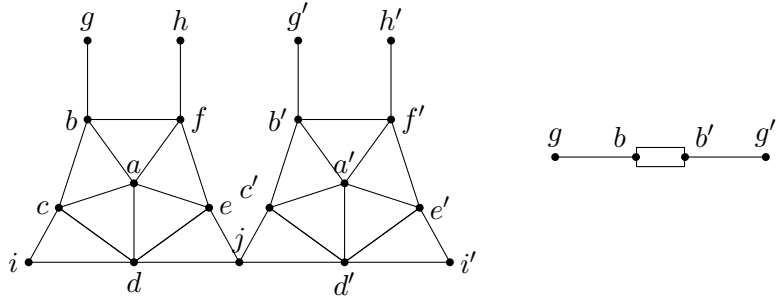
\begin{figure}[htb]
  \begin{center}
    \begin{tikzpicture}[line cap=round,line join=round,x=.7cm,y=.7cm]
      \coordinate (a) at (-2,1.5);
      \coordinate (b) at ($ (a) + (126:1.5) $);
      \coordinate (c) at ($ (a) + (-162:1.5) $);
      \coordinate (d) at ($ (a) + (-90:1.5) $);
      \coordinate (e) at ($ (a) + (-18:1.5) $);
      \coordinate (f) at ($ (a) + (54:1.5) $);
      \coordinate (g) at ($ (b) + (90:1.5) $);
      \coordinate (h) at ($ (f) + (90:1.5) $);
      \coordinate (i) at ($ (d) + (-2,0) $);
      \coordinate (j) at ($ (d) + (2,0) $);
      \coordinate (a') at (2,1.5);
      \coordinate (b') at ($ (a') + (126:1.5) $);
      \coordinate (c') at ($ (a') + (-162:1.5) $);
      \coordinate (d') at ($ (a') + (-90:1.5) $);
      \coordinate (e') at ($ (a') + (-18:1.5) $);
      \coordinate (f') at ($ (a') + (54:1.5) $);
      \coordinate (g') at ($ (b') + (90:1.5) $);
      \coordinate (h') at ($ (f') + (90:1.5) $);
      \coordinate (i') at ($ (d') + (2,0) $);

      \fill (a) circle (1.5pt);
      \fill (b) circle (1.5pt);
      \fill (c) circle (1.5pt);
      \fill (d) circle (1.5pt);
      \fill (e) circle (1.5pt);
      \fill (f) circle (1.5pt);
      \fill (g) circle (1.5pt);
      \fill (h) circle (1.5pt);
      \fill (i) circle (1.5pt);
      \fill (j) circle (1.5pt);
       \fill (a') circle (1.5pt);
      \fill (b') circle (1.5pt);
      \fill (c') circle (1.5pt);
      \fill (d') circle (1.5pt);
      \fill (e') circle (1.5pt);
      \fill (f') circle (1.5pt);
      \fill (g') circle (1.5pt);
      \fill (h') circle (1.5pt);
      \fill (i') circle (1.5pt);

      \draw (b) -- (c) -- (d) -- (e) -- (f) -- (b);
      \draw (c) -- (i) -- (d) -- (c);
      \draw (d) -- (j) -- (e) -- (d);
      
       \draw (b') -- (c') -- (d') -- (e') -- (f') -- (b');
      \draw (c') -- (j) -- (d') -- (c');
      \draw (d') -- (i') -- (e') -- (d');

      \draw (a) -- (b);
      \draw (a) -- (c);
      \draw (a) -- (d);
      \draw (a) -- (e);
      \draw (a) -- (f);
      \draw (b) -- (g);
      \draw (f) -- (h);

      \draw (a') -- (b');
      \draw (a') -- (c');
      \draw (a') -- (d');
      \draw (a') -- (e');
      \draw (a') -- (f');
      \draw (b') -- (g');
      \draw (f') -- (h');
      
      \draw (a) node[above] {$a$};
      \draw (b) node[left] {$b$};
      \draw (c) node[left] {$c$};
      \draw (d) node[below] {$d$};
      \draw (e) node[right] {$e$};
      \draw (f) node[right] {$f$}; 
      \draw (g) node[above] {$g$};
      \draw (h) node[above] {$h$}; 
      \draw (i) node[left] {$i$};
      \draw (j) node[above] {$j$}; 
      \draw (a') node[above] {$a'$};
      \draw (b') node[left] {$b'$};
      \draw (c') node[above left] {$c'$};
      \draw (d') node[below] {$d'$};
      \draw (e') node[right] {$e'$};
      \draw (f') node[right] {$f'$}; 
      \draw (g') node[above] {$g'$};
      \draw (h') node[above] {$h'$}; 
      \draw (i') node[right] {$i'$};

      \node[draw] (inv) at (8,2) {$~~~$};
      \coordinate (gg) at (6,2);
      \coordinate (gg') at (10,2);
      \fill (gg) circle (1.5pt);
      \fill (gg') circle (1.5pt);
       \fill (inv.east) circle (1.5pt);
       \fill (inv.west) circle (1.5pt);
       \draw (gg) node[above] {$g$};
       \draw (gg') node[above] {$g'$};
      \draw (inv.east) node[above right] {$b'$};
      \draw (inv.west) node[above left] {$b$};
      
      \draw (gg) -- (inv.west);
      \draw (inv.east) -- (gg');      

    \end{tikzpicture}
  \end{center}
  \caption{The graph \textsc{Inv} and its abbreviation.}
  \label{fig:inv}
\end{figure}

  Let us call \textsc{Inv} (\textsc{Inv} stands for inverter) the
  graph obtained from $I$ and $I'$ by identifying the nodes $j$ and
  $j'$. We use $j$ for the name of the resulting node, see Figure
  \ref{fig:inv}.
  
From the discussion above there are only two directed graphs that for
convenience we call $\overleftrightarrow{\textsc{Inv}}$ and
$\overline{\textsc{Inv}}$, such that
$\textsc{Inv}=I(\overleftrightarrow{\textsc{Inv}})=I(\overline{\textsc{Inv}})$. The
graph $\overleftrightarrow{\textsc{Inv}}$
(respectively $\overline{\textsc{Inv}}$) is obtained from $D_1$ and
$D'_1$ (respectively $D_2$ and $D'_2$) by identifying $j$ and $j'$. Notice
that the other possibilities of identifying $j$ and $j'$ will not lead
to the graph $\textsc{Inv}$.
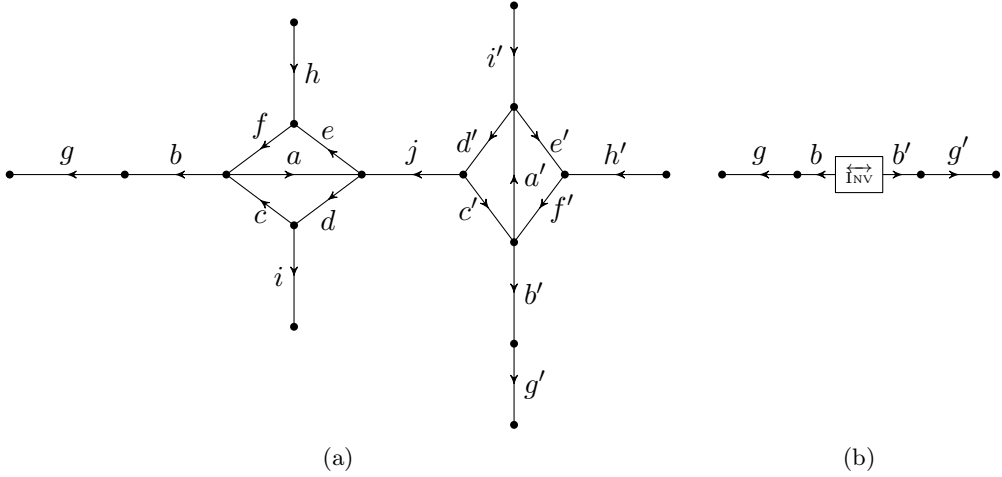
\begin{figure}[htb]
  \begin{center}
    \subfigure[]{\begin{tikzpicture}[line cap=round,line join=round,x=.9cm,y=.9cm]
      \coordinate (x) at (4,2.25);
      \coordinate (y) at (3,3);
      \coordinate (z) at (2,2.25);
      \coordinate (t) at (3,1.5);
      \coordinate (u) at (3,0);
      
      \coordinate (p) at (3,4.5);
      \coordinate (q) at (0.5,2.25);
      \coordinate (w) at (5.5,2.25);
      \coordinate (a) at (6.25,4.75);
      \coordinate (b) at (6.25,3.25);
      \coordinate (c) at (6.25,1.25);
      \coordinate (d) at (6.25,-0.25);
      \coordinate (e) at (7,2.25);
      \coordinate (f) at (8.5,2.25);
      \coordinate (v) at (-1.2,2.25);
      \coordinate (g) at (6.25,-1.45);
      
      \begin{scope}[decoration={
          markings,
          mark=at position 0.5 with {\arrow{stealth'}}}
        ] 
        
        \draw[postaction={decorate}] (t) -- (z) node[midway,left, below]{$c$};
        \draw[postaction={decorate}] (x) -- (t) node[midway,right, below]{$d$};
        \draw[postaction={decorate}] (z) -- (q) node[midway,above]{$b$};
        \draw[postaction={decorate}] (x) -- (y) node[midway,right, above]{$e$};
        \draw[postaction={decorate}] (z) -- (x) node[midway,above]{$a$};
        \draw[postaction={decorate}] (y) -- (z) node[midway,left, above]{$f$};
        \draw[postaction={decorate}] (t) -- (u) node[midway,left]{$i$};
        \draw[postaction={decorate}] (w) -- (x) node[midway, above]{$j$};
        \draw[postaction={decorate}] (q) -- (v) node[midway,above]{$g$};
        \draw[postaction={decorate}] (p) -- (y) node[midway,right]{$h$};
        
        \draw[postaction={decorate}] (c) -- (b) node[midway, right]{$a'$};
        \draw[postaction={decorate}] (b) -- (w) node[midway,left]{$d'$};
        \draw[postaction={decorate}] (w) -- (c) node[midway,left]{$c'$};
        \draw[postaction={decorate}] (d) -- (g) node[midway,right]{$g'$};
        \draw[postaction={decorate}] (c) -- (d) node[midway,right]{$b'$};
        \draw[postaction={decorate}] (b) -- (e) node[midway,right]{$e'$};
        \draw[postaction={decorate}] (a) -- (b) node[midway, left]{$i'$};
        \draw[postaction={decorate}] (e) -- (c) node[midway, right]{$f'$};
        \draw[postaction={decorate}] (f) -- (e) node[midway,above]{$h'$};
        
      \end{scope}
      
      \fill (x) circle (1.5pt);
      \fill (y) circle (1.5pt);
      \fill (z) circle (1.5pt);
      \fill (t) circle (1.5pt);
      \fill (u) circle (1.5pt);
      \fill (p) circle (1.5pt);
      \fill (q) circle (1.5pt);
      \fill (v) circle (1.5pt);
      \fill (w) circle (1.5pt);
      \fill (a) circle (1.5pt);
      \fill (b) circle (1.5pt);
      \fill (c) circle (1.5pt);
      \fill (d) circle (1.5pt);
      \fill (e) circle (1.5pt);
      \fill (f) circle (1.5pt);
      \fill (g) circle (1.5pt);
      \end{tikzpicture} \label{fig:inv_graph1}
    } \quad \subfigure[]{\begin{tikzpicture}[line cap=round,line join=round,x=.5cm,y=.5cm]
      \node[draw] (inv) at (12,3) {\tiny $\overleftrightarrow{\textsc{Inv}}$};
      \coordinate (a') at ($ (inv.east) + (1,0) $);
      \coordinate (b') at ($ (inv.east) + (3,0) $);
      \coordinate (a'') at ($ (inv.west) + (-1,0) $);
      \coordinate (b'') at ($ (inv.west) + (-3,0) $);
      \coordinate (cadrage) at ($ (inv) + (0,-6.7) $);
      \fill (a'') circle (1.5pt);
      \fill (b'') circle (1.5pt);
      \fill (a') circle (1.5pt);
      \fill (b') circle (1.5pt);
      \fill[color=white] (cadrage) circle (1pt);
        
      \begin{scope}[decoration={
          markings,
          mark=at position 0.5 with {\arrow{stealth'}}}
        ] 
        \draw[postaction={decorate}] (inv.east) -- (a') node[midway, above]{$b'$};
        \draw[postaction={decorate}] (a') -- (b') node[midway, above]{$g'$};
        \draw[postaction={decorate}] (inv.west) -- (a'') node[midway, above]{$b$};
        \draw[postaction={decorate}] (a'') -> (b'') node[midway, above]{$g$};
      \end{scope}
      
      \end{tikzpicture} \label{fig:inv_abbrev1}
    }
    \end{center}
    \caption{The graph $\protect\overleftrightarrow{\textsc{Inv}}$ (a) and its abbreviation (b).}
    \label{fig:inv_depic1}
  \end{figure}
  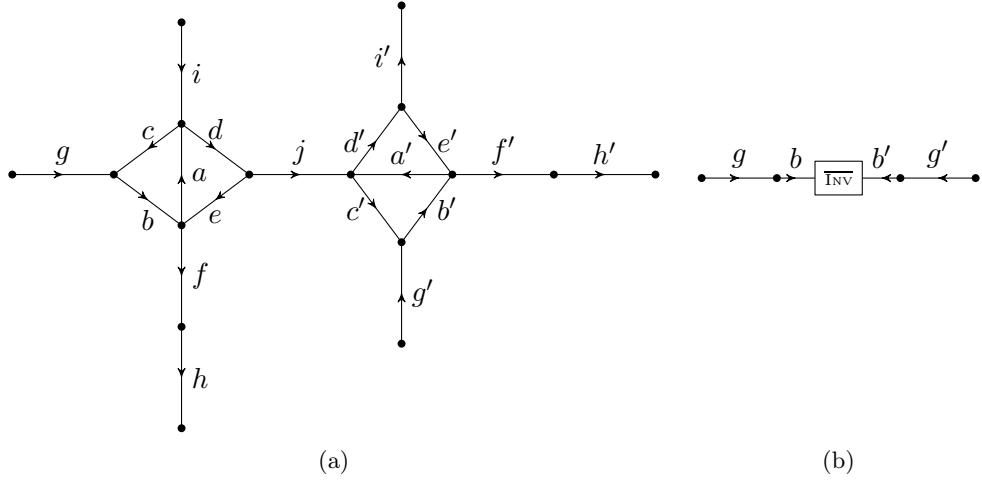
\begin{figure}[htb]
    \begin{center}
     \subfigure[]{\begin{tikzpicture}[line cap=round,line join=round,x=.9cm,y=.9cm]
        \coordinate (x) at (4,2.25);
        \coordinate (y) at (3,3);
        \coordinate (z) at (2,2.25);
        \coordinate (t) at (3,1.5);
        \coordinate (u) at (3,0);
        \coordinate (v) at (3,-1.5);
        \coordinate (p) at (3,4.5);
        \coordinate (q) at (0.5,2.25);
        \coordinate (w) at (5.5,2.25);
        \coordinate (a) at (6.25,4.75);
        \coordinate (b) at (6.25,3.25);
        \coordinate (c) at (6.25,1.25);
        \coordinate (d) at (6.25,-0.25);
        \coordinate (e) at (7,2.25);
        \coordinate (f) at (8.5,2.25);
        \coordinate (g) at (10,2.25);

        \begin{scope}[decoration={
            markings,
            mark=at position 0.5 with {\arrow{stealth'}}}
          ] 
        
        \draw[postaction={decorate}] (y) -- (z) node[midway,left, above]{$c$};
        \draw[postaction={decorate}] (y) -- (x) node[midway,right, above]{$d$};
        \draw[postaction={decorate}] (z) -- (t) node[midway,left, below]{$b$};
        \draw[postaction={decorate}] (x) -- (t) node[midway,right, below]{$e$};
        \draw[postaction={decorate}] (t) -- (y) node[midway,right]{$a$};
        \draw[postaction={decorate}] (t) -- (u) node[midway,right]{$f$};
        \draw[postaction={decorate}] (p) -- (y) node[midway,right]{$i$};
        \draw[postaction={decorate}] (x) -- (w) node[midway, above]{$j$};
        \draw[postaction={decorate}] (q) -- (z) node[midway,above]{$g$};
        \draw[postaction={decorate}] (u) -- (v) node[midway,right]{$h$};
        
        \draw[postaction={decorate}] (e) -- (w) node[midway, above]{$a'$};
        \draw[postaction={decorate}] (w) -- (b) node[midway,left]{$d'$};
        \draw[postaction={decorate}] (w) -- (c) node[midway,left]{$c'$};
        \draw[postaction={decorate}] (d) -- (c) node[midway,right]{$g'$};
        \draw[postaction={decorate}] (c) -- (e) node[midway,right]{$b'$};
        \draw[postaction={decorate}] (b) -- (e) node[midway,right]{$e'$};
        \draw[postaction={decorate}] (b) -- (a) node[midway, left]{$i'$};
        \draw[postaction={decorate}] (e) -- (f) node[midway,above]{$f'$};
        \draw[postaction={decorate}] (f) -- (g) node[midway,above]{$h'$};

      \end{scope}
        
        \fill (x) circle (1.5pt);
        \fill (y) circle (1.5pt);
        \fill (z) circle (1.5pt);
        \fill (t) circle (1.5pt);
        \fill (u) circle (1.5pt);
        \fill (p) circle (1.5pt);
        \fill (q) circle (1.5pt);
        \fill (v) circle (1.5pt);
        \fill (w) circle (1.5pt);
        \fill (a) circle (1.5pt);
        \fill (b) circle (1.5pt);
        \fill (c) circle (1.5pt);
        \fill (d) circle (1.5pt);
        \fill (e) circle (1.5pt);
        \fill (f) circle (1.5pt);
        \fill (g) circle (1.5pt);
      \end{tikzpicture}
    }\quad \subfigure[]{\begin{tikzpicture}[line cap=round,line join=round,x=.5cm,y=.5cm]
      \node[draw] (inv) at (12,3) {\tiny $\overline{\textsc{Inv}}$};
      \coordinate (a') at ($ (inv.east) + (1,0) $);
      \coordinate (b') at ($ (inv.east) + (3,0) $);
      \coordinate (a'') at ($ (inv.west) + (-1,0) $);
      \coordinate (b'') at ($ (inv.west) + (-3,0) $);
      \coordinate (cadrage) at ($ (inv) + (0,-6.7) $);
      \fill (a'') circle (1.5pt);
      \fill (b'') circle (1.5pt);
      \fill (a') circle (1.5pt);
      \fill (b') circle (1.5pt);
      \fill[color=white] (cadrage) circle (1pt);
        
      \begin{scope}[decoration={
          markings,
          mark=at position 0.5 with {\arrow{stealth'}}}
        ] 
        \draw[postaction={decorate}] (a') -- (inv.east)  node[midway, above]{$b'$};
        \draw[postaction={decorate}] (b') -- (a') node[midway, above]{$g'$};
        \draw[postaction={decorate}] (a'') -- (inv.west) node[midway, above]{$b$};
        \draw[postaction={decorate}] (b'') -> (a'') node[midway, above]{$g$};
      \end{scope}
      
      \end{tikzpicture}
    }
    \end{center}
 \caption{The graph $\protect\overline{\textsc{Inv}}$ (a) and its abbreviation (b).}
    \label{fig:inv_depic2}
  \end{figure}

\begin{figure}[htb]
  \begin{center}
    \begin{tikzpicture}[line cap=round,line join=round,x=.65cm,y=.65cm]
      \coordinate (a) at (-2,1.5);
      \coordinate (b) at ($ (a) + (126:1.5) $);
      \coordinate (c) at ($ (a) + (-162:1.5) $);
      \coordinate (d) at ($ (a) + (-90:1.5) $);
      \coordinate (e) at ($ (a) + (-18:1.5) $);
      \coordinate (f) at ($ (a) + (54:1.5) $);
      \coordinate (g) at ($ (b) + (90:1.5) $);
      \coordinate (h) at ($ (f) + (90:1.5) $);
      \coordinate (i) at ($ (d) + (-2,0) $);
      \coordinate (j) at ($ (d) + (2,0) $);
      \coordinate (a') at (2,1.5);
      \coordinate (b') at ($ (a') + (126:1.5) $);
      \coordinate (c') at ($ (a') + (-162:1.5) $);
      \coordinate (d') at ($ (a') + (-90:1.5) $);
      \coordinate (e') at ($ (a') + (-18:1.5) $);
      \coordinate (f') at ($ (a') + (54:1.5) $);
      \coordinate (g') at ($ (b') + (90:1.5) $);
      \coordinate (h') at ($ (f') + (90:1.5) $);
      \coordinate (i') at ($ (d') + (2,0) $);
      \coordinate (am) at (10,1.5);
      \coordinate (bm) at ($ (am) + (126:1.5) $);
      \coordinate (cm) at ($ (am) + (-162:1.5) $);
      \coordinate (dm) at ($ (am) + (-90:1.5) $);
      \coordinate (em) at ($ (am) + (-18:1.5) $);
      \coordinate (fm) at ($ (am) + (54:1.5) $);
      \coordinate (gm) at ($ (bm) + (90:1.5) $);
      \coordinate (hm) at ($ (fm) + (90:1.5) $);
      \coordinate (im) at ($ (dm) + (-2,0) $);
      \coordinate (im1) at ($ (dm) + (2,0) $);

      \fill (a) circle (1.5pt);
      \fill (b) circle (1.5pt);
      \fill (c) circle (1.5pt);
      \fill (d) circle (1.5pt);
      \fill (e) circle (1.5pt);
      \fill (f) circle (1.5pt);
      \fill (g) circle (1.5pt);
      \fill (h) circle (1.5pt);
      \fill (i) circle (1.5pt);
      \fill (j) circle (1.5pt);
       \fill (a') circle (1.5pt);
      \fill (b') circle (1.5pt);
      \fill (c') circle (1.5pt);
      \fill (d') circle (1.5pt);
      \fill (e') circle (1.5pt);
      \fill (f') circle (1.5pt);
      \fill (g') circle (1.5pt);
      \fill (h') circle (1.5pt);
      \fill (i') circle (1.5pt);
      \fill (am) circle (1.5pt);
      \fill (bm) circle (1.5pt);
      \fill (cm) circle (1.5pt);
      \fill (dm) circle (1.5pt);
      \fill (em) circle (1.5pt);
      \fill (fm) circle (1.5pt);
      \fill (gm) circle (1.5pt);
      \fill (hm) circle (1.5pt);
      \fill (im) circle (1.5pt);
      \fill (im1) circle (1.5pt);

      \draw (b) -- (c) -- (d) -- (e) -- (f) -- (b);
      \draw (c) -- (i) -- (d) -- (c);
      \draw (d) -- (j) -- (e) -- (d);
      
       \draw (b') -- (c') -- (d') -- (e') -- (f') -- (b');
      \draw (c') -- (j) -- (d') -- (c');
      \draw (d') -- (i') -- (e') -- (d');

      \draw (bm) -- (cm) -- (dm) -- (em) -- (fm) -- (bm);
      \draw (cm) -- (im) -- (dm) -- (cm);
      \draw (dm) -- (im1) -- (em) -- (dm);

      \draw (a) -- (b);
      \draw (a) -- (c);
      \draw (a) -- (d);
      \draw (a) -- (e);
      \draw (a) -- (f);
      \draw (b) -- (g);
      \draw (f) -- (h);
      
      \draw (am) -- (bm);
      \draw (am) -- (cm);
      \draw (am) -- (dm);
      \draw (am) -- (em);
      \draw (am) -- (fm);
      \draw (bm) -- (gm);
      \draw (fm) -- (hm);

      \draw (a') -- (b');
      \draw (a') -- (c');
      \draw (a') -- (d');
      \draw (a') -- (e');
      \draw (a') -- (f');
      \draw (b') -- (g');
      \draw (f') -- (h');

      \draw[dotted] (i') -- (im);
      
      \begin{scriptsize}
      \draw ($ (a) + (25:0.4) $) node {$a_1^i$};
      \draw (b) node[left] {$b_1^i$};
      \draw (c) node[left] {$c_1^i$};
      \draw (d) node[below] {$d_1^i$};
      \draw (e) node[right] {$e_1^i$};
      \draw (f) node[right] {$f_1^i$}; 
      \draw (g) node[above] {$g_1^i$};
      \draw (h) node[above] {$h_1^i$}; 
      \draw (i) node[left] {$i_1^i$};
      \draw (j) node[below] {$i_2^i$}; 
      \draw ($ (a') + (25:0.4) $) node {$a_2^i$};
      \draw (b') node[left] {$b_2^i$};
      \draw (c') node[above left] {$c_2^i$};
      \draw (d') node[below] {$d_2^i$};
      \draw (e') node[right] {$e_2^i$};
      \draw (f') node[right] {$f_2^i$}; 
      \draw (g') node[above] {$g_2^i$};
      \draw (h') node[above] {$h_2^i$}; 
      \draw (i') node[below] {$i_3^i$};
      \draw ($ (am) + (20:0.5) $) node {$a_m^i$};
      \draw (bm) node[left] {$b_m^i$};
      \draw (cm) node[above left] {$c_m^i$};
      \draw (dm) node[below] {$d_m^i$};
      \draw (em) node[right] {$e_m^i$};
      \draw (fm) node[right] {$f_m^i$}; 
      \draw (gm) node[above] {$g_m^i$};
      \draw (hm) node[above] {$h_m^i$}; 
      \draw (im) node[below] {$i_m^i$};
      \draw (im1) node[right] {$i_{m+1}^i$};
      \end{scriptsize}

    \end{tikzpicture}
  \end{center}
  \caption{Graph for every variable $x_i$, $\textsc{Gad}^1_i$.}
  \label{fig:gadget1}
\end{figure}
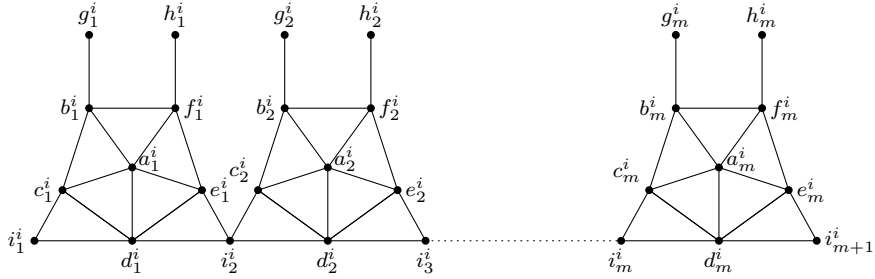

  Now we are ready to build graphs $\textsc{Gad}^1_i$ that correspond
  to the variables $x_i$, for each $i$ in $\{1,\dots,n\}$.  For each
  variable $x_i$ we construct $m$ copies of the graph $I$, where the
  nodes $a,\ldots,j$ of each copy are renamed, respectively,
  $a^i_1,\ldots,j^i_1$ up to $a^i_m,\ldots,j^i_m$. The graph
  $\textsc{Gad}^1_i$ is obtained by identifying the node $j^i_l$ with
  $i^i_{l+1}$ and we call $i^i_{l+1}$ the resulting node, for
  $l=1,\dots,m-1$. Also we rename the node $j^i_m$ by $i^i_{m+1}$, see
  Figure \ref{fig:gadget1}.

The discussion above implies the following lemma, see Figures
\ref{fig:inv_depic1} and \ref{fig:inv_depic2}.
 
\begin{lemma}
  \label{GAD1}
  For each directed graph $D$ with $I(D)=\textsc{Gad}^1_i$ exactly one
  of the following two assumptions holds:
  
  \begin{itemize}
  \item[(i)] $h(b^i_j)=t(g^i_j)$ and $t(f^i_j)=h(h^i_j)$ for each
    $j=1,\dots,m$,
  \item[(ii)] $t(b^i_j)=h(g^i_j)$ and $h(f^i_j)=t(h^i_j)$ for each
    $j=1,\dots,m$.
  \end{itemize}
  
\end{lemma}

\subsection{The construction of the graphs $\textsc{Gad}^2_j$}

\begin{figure}[htb]
  \begin{center}
    \begin{tikzpicture}[line cap=round,line join=round,x=.5cm,y=.5cm]
      \coordinate (r') at (0,0);
      \coordinate (r) at ($ (r') + (0,-2) $);
      \coordinate (a) at ($ (r) + (-120:2) $);
      \coordinate (f) at ($ (r) + (-60:2) $);
      \node[draw,rotate=60] (g) at ($ (a) + (-120:3) $) {$~~~$};
      \node[draw,rotate=-60] (h) at ($ (f) + (-60:3) $) {$~~~$};
      \coordinate (b) at ($ (a) + (-120:6) $);
      \coordinate (e) at ($ (f) + (-60:6) $);
      \coordinate (s) at ($ (b) + (-120:2) $);
      \coordinate (c) at ($ (b) + (-60:2) $);
      \coordinate (d) at ($ (e) + (-120:2) $);
      \coordinate (t) at ($ (e) + (-60:2) $);
      \coordinate (s') at ($ (s) + (-120:2) $);
      \coordinate (t') at ($ (t) + (-60:2) $);

      \coordinate (cen) at (12,-7);
      \coordinate (n) at ($ (cen) + (90:2) $);
      \coordinate (sw) at ($ (cen) + (-150:2) $);
      \coordinate (se) at ($ (cen) + (-30:2) $);
      \coordinate (n') at ($ (cen) + (90:4) $);
      \coordinate (sw') at ($ (cen) + (-150:4) $);
      \coordinate (se') at ($ (cen) + (-30:4) $);


      \draw (r) -- (r');
      \draw (t) -- (t');
      \draw (s) -- (s');
      \draw (r) -- (a) -- (f) -- (r);
      \draw (s) -- (b) -- (c) -- (s);
      \draw (t) -- (d) -- (e) -- (t);
      \draw (d) -- (c);
      \draw (g.east) -- (a);
      \draw (g.west) -- (b);
      \draw (h.east) -- (e);
      \draw (h.west) -- (f);
      
      \fill (a) circle (1.5pt);
      \fill (b) circle (1.5pt);
      \fill (c) circle (1.5pt);
      \fill (d) circle (1.5pt);
      \fill (e) circle (1.5pt);
      \fill (f) circle (1.5pt);
      \fill (r) circle (1.5pt);
      \fill (r') circle (1.5pt);
      \fill (s) circle (1.5pt);
      \fill (s') circle (1.5pt);
      \fill (t) circle (1.5pt);
      \fill (t') circle (1.5pt);
      \fill (g.east) circle (1.5pt);
      \fill (h.east) circle (1.5pt);
      \fill (g.west) circle (1.5pt);
      \fill (h.west) circle (1.5pt);
      
      \draw (g.east) node[above] {$a'$};
      \draw (g.west) node[left] {$b'$};
      \draw (h.east) node[right] {$e'$};
      \draw (h.west) node[above] {$f'$};
      \draw (a) node[above left] {$a$};
      \draw (b) node[above left] {$b$};
      \draw (c) node[above right] {$c$};
      \draw (d) node[above left] {$d$};
      \draw (e) node[above right] {$e$};
      \draw (f) node[above right] {$f$};
      \draw (r) node[above right] {$r$};
      \draw (r') node[above] {$r'$};
      \draw (s) node[above left] {$s$};
      \draw (s') node[below left] {$s'$};
      \draw (t) node[right] {$t$};
      \draw (t') node[below right] {$t'$};
    \end{tikzpicture}
  \end{center}
  \caption{The graph $Gad^2_j$ corresponding to the clause $C_j$.}
  \label{fig:gadget2}
\end{figure}
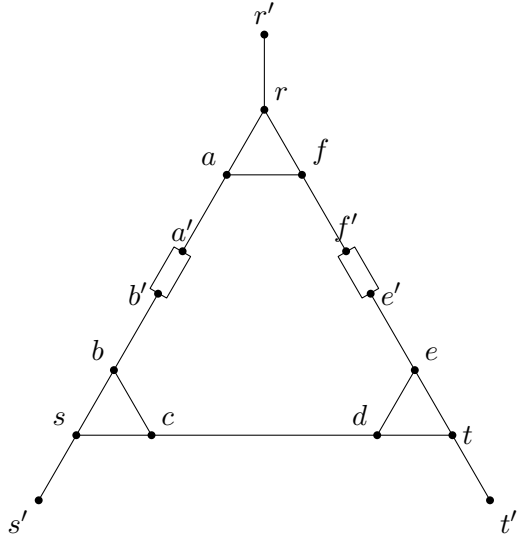

We will use the graph $\textsc{Inv}$ of the previous subsection to
construct $\textsc{Gad}^2_j$.  We have three triangles
$\Delta_1=\{r_j,a_j,f_j\}$, $\Delta_2=\{s_j,b_j,c_j\}$ and
$\Delta_3=\{t_j,e_j,d_j\}$ with the addition of a branch pending from
each node of these triangles, that is we add the edges $r_jr'_j$,
$a_jf'_j$; $s_js'_j$, $b_jb'_j$; $t_jt'_j$, $e_je'_j$, $d_jd'_j$.
These triangles are connected using their branches. We choose one
triangle say $\Delta_1$ and we connect it to $\Delta_2$ and $\Delta_3$
via two graphs identical to $\textsc{Inv}$ using two of its branches
$a_ja'_j$ and $f_jf'_j$. The triangles $\Delta_2$ and $\Delta_3$ are
connected by identifying the branches $c_jc'_j$ and $d_jd'_j$ (the
nodes $c'_j$ and $d'_j$ are removed), see Figure \ref{fig:gadget2}.

Before establishing the main lemma of this section let us notice the
following remark.

\begin{figure}[htb]
  \begin{center}
    \begin{tikzpicture}[line cap=round,line join=round,x=.6cm,y=.6cm]
      \coordinate (a) at (0,2);
      \coordinate (b) at ($ (a) + (-1.5,-2) $);
      \coordinate (c) at ($ (a) + (0,-2) $);
      \coordinate (d) at ($ (a) + (1.5,-2) $);
      \coordinate (e) at (4,1.5);
      \coordinate (f) at ($ (e) + (0,1.5) $);
      \coordinate (g) at ($ (e) + (-0.75,-1.5) $);
      \coordinate (h) at ($ (e) + (0.75,-1.5) $);
      \coordinate (k) at (8,2);
      \coordinate (i) at ($ (k) + (-1.3,-2) $);
      \coordinate (j) at ($ (k) + (1.3,-2) $);
      \coordinate (m) at (12,2);
      \coordinate (l) at ($ (m) + (0,2) $);
      \coordinate (n) at ($ (m) + (0,-2) $);
      \coordinate (o) at ($ (b) + (0,-2) $);
      \coordinate (p) at ($ (c) + (0,-2) $);
      \coordinate (q) at ($ (d) + (0,-2) $);
      \coordinate (r) at ($ (g) + (0,-2) $);
      \coordinate (s) at ($ (h) + (0,-2) $);
      \coordinate (t) at ($ (f) + (0,2) $);
      \coordinate (u) at ($ (i) + (-120:2) $);
      \coordinate (w) at ($ (j) + (-60:2) $);
      \coordinate (v) at ($ (k) + (90:2) $);
       \coordinate (p1) at ($ (c) + (0,-2.5) $);
       \coordinate (s1) at ($ (h) + (-0.5,-2.5) $);
        \coordinate (s1) at ($ (h) + (-0.5,-2.5) $);
       \coordinate (v1) at (8,-2.5);

      \fill (a) circle (1.5pt);
      \fill (b) circle (1.5pt);
      \fill (c) circle (1.5pt);
      \fill (d) circle (1.5pt);
      \fill (e) circle (1.5pt);
      \fill (f) circle (1.5pt);
      \fill (g) circle (1.5pt);
      \fill (h) circle (1.5pt);
      \fill (i) circle (1.5pt);
      \fill (j) circle (1.5pt);
      \fill (k) circle (1.5pt);
       \fill (o) circle (1.5pt);
      \fill (p) circle (1.5pt);
      \fill (q) circle (1.5pt);
      \fill (r) circle (1.5pt);
      \fill (s) circle (1.5pt);
      \fill (t) circle (1.5pt);
      \fill (u) circle (1.5pt);
      \fill (v) circle (1.5pt);
      \fill (w) circle (1.5pt);
      
      \begin{scope}[decoration={
            markings,
            mark=at position 0.5 with {\arrow{stealth'}}}
          ] 

      \draw[postaction={decorate}] (a) -- (b);
      \draw[postaction={decorate}] (a) -- (c);
      \draw[postaction={decorate}] (a) -- (d);

      \draw[postaction={decorate},dashed] (b) -- (o);
      \draw[postaction={decorate},dashed] (c) -- (p);
      \draw[postaction={decorate},dashed] (d) -- (q);
       \draw (p1) node[below] {$D_1$};
      
      \draw[postaction={decorate}] (f) -- (e);
      \draw[postaction={decorate}] (e) -- (h);
      \draw[postaction={decorate}] (e) -- (g);

      \draw[postaction={decorate},dashed] (t) -- (f);
      \draw[postaction={decorate},dashed] (g) -- (r);
      \draw[postaction={decorate},dashed] (h) -- (s);
       \draw (s1) node[below] {$D_2$};

      \draw[postaction={decorate}] (i) -- (j);
      \draw[postaction={decorate}] (j) -- (k);
      \draw[postaction={decorate}] (k) -- (i);
      
      \draw[postaction={decorate},dashed] (u) -- (i);
      \draw[postaction={decorate},dashed] (w) -- (j);
      \draw[postaction={decorate},dashed] (v) -- (k);
       \draw (v1) node[below] {$D_3$};

    \end{scope}
      
    \end{tikzpicture}
  \end{center}
  \caption{The solid line are the nodes of the  triangle $\Delta$ and the dashed lines are its pending nodes.}
  \label{fig:k3s}
\end{figure}
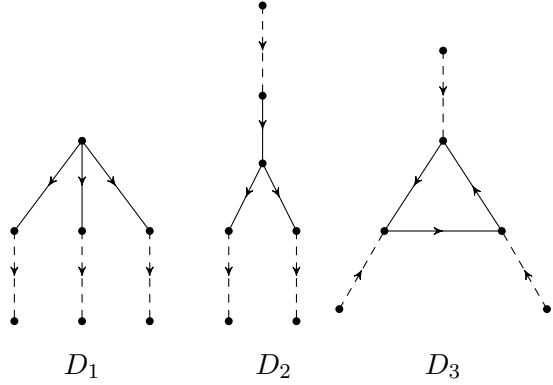

\begin{obs}
  \label{Delta}
  Call $\Delta$ the undirected graph defined by a triangle with three
  branches pending from each of its three nodes.  There are only three
  possible directed graphs, $D_1$, $D_2$, $D_3$ such that
  $\Delta=I(D_1)=I(D_2)=I(D_3)$, as shown in Figure \ref{fig:k3s}
\end{obs}

\begin{lemma}
  \label{GAD2}
  Let $D$ be a directed graph such that $I(D)=\textsc{Gad}^2_j$. Then
  the two following assumptions hold:
  
  \begin{itemize}
  \item[(i)] The arcs $r'_j$, $s'_j$ and $t'_j$ cannot all enter the
    arcs $r_j$, $s_j$ and $t_j$, respectively,
  \item[(ii)] Any other configuration for theses three adjacencies is
    possible.
  \end{itemize}
  
\end{lemma}
\proof (i) Assume that in $D$ all the arcs $r'_j$, $s'_j$ and $t'_j$
enter the arcs $r_j$, $s_j$ and $t_j$, respectively. That is
$h(r'_j)=t(r_j)$, $h(s'_j)=t(s_j)$ and $h(t'_j)=t(t_j)$.  Therefore,
from Remark \ref{Delta}, none of the triangles $\Delta_1$, $\Delta_2$
and $\Delta_3$ has the configuration $D_1$ of Figure \ref{fig:k3s}.
We may also check that $\Delta_2$ and $\Delta_3$ cannot have the same
configuration, that is cannot be both of the form $D_2$ or $D_3$ of
Figure \ref{fig:k3s}.  Hence assume that $\Delta_2$ has the form of
$D_2$ and $\Delta_3$ has the form of $D_3$.  Since the arc $e_j'$ must
enter the arc $e_j$ it follows from Lemma \ref{GAD1} that
 \begin{eqnarray}
 \label{1r}
 f'_j \mbox{ must enter the arc } f_j.
 \end{eqnarray}
 Since $c_j$ enters the arc $d_j$ and using the fact that $s'_j$
 enters $s_j$ we conclude that $b_j$ must enter the arc $b'_j$. Now
 using again Lemma \ref{GAD1} we obtain
 \begin{eqnarray}
 \label{2r}
 a_j  \mbox{ must enter the arc } a'_j.
 \end{eqnarray}

 Now combining the two facts \eqref{1r} and \eqref{2r} we have that
 $\Delta_1$ must be of the form $D_2$ of Figure \ref{fig:k3s} and that
 the arc $r_j$ must enter $r'_j$, which is not possible.
 
 (ii) The proof of this assumption is presented in the appendix. It
 lists all the possible configurations of $D$ in Figure \ref{fig:preimages_gad2}

 \endproof

 \subsection{The construction of the graph $G_F$}

 Let $F=C_1\wedge\dots\wedge C_m$, where each clause
 $C_j=\lambda_{j_1}\vee\lambda_{j_2}\vee\lambda_{j_3}$, for
 $j$ in $\{1,\dots,m\}$. Each $\lambda_{j_k}$ correspond to the variable
 $x_{j_k}$ or its negation $\bar{x}_{j_k}$.  From $F$ we construct an
 undirected graph $G_F$ as follows. Let $\textsc{Gad}^1_i$ be the
 undirected graph associated with each Boolean variable $x_i$,
 $i=1,\dots,n$. And let $\textsc{Gad}^2_j$ be the undirected graph
 associated with each clause $C_j$, $j=1,\dots,m$.  In the
 construction of $G_F$, there is no connection between the graphs
 $\textsc{Gad}^1_i$ themselves and between the graphs
 $\textsc{Gad}^2_j$ themselves. The only connections are between
 $\textsc{Gad}^1_i$ and $\textsc{Gad}^2_j$ where the variable $x_i$ or
 $\bar{x}_i$ appears in the clause $C_j$. Moreover these connections
 are made through their branches. Specifically, for each clause
 $C_j=\lambda_{j_1}\vee\lambda_{j_2}\vee\lambda_{j_3}$ we do the
 following:

\begin{itemize}

\item if $\lambda_{j_1} = x_{j_1}$, we identify vertex $r_j$ with vertex
  $g_j^{j_1}$ and vertex $r'_j$ with vertex $b_j^{j_1}$,
\item if $\lambda_{j_1} = \bar{x}_{j_1}$, we identify vertex $r_j$ with vertex
  $h_j^{j_1}$ and vertex $r'_j$ with vertex $f_j^{j_1}$,
\item if $\lambda_{j_2} = x_{j_2}$, we identify vertex $s_j$ with vertex
  $g_j^{j_2}$ and vertex $s'_j$ with vertex $b_j^{j_2}$,
\item if $\lambda_{j_2} = \bar{x}_{j_2}$, we identify vertex $s_j$ with vertex
  $h_j^{j_2}$ and vertex $s'_j$ with vertex $f_j^{j_2}$,
\item if $\lambda_{j_3} = x_{j_3}$, we identify vertex $t_j$ with vertex
  $g_j^{j_3}$ and vertex $t'_j$ with vertex $b_j^{j_3}$,
\item if $\lambda_{j_3} = \bar{x}_{j_3}$, we identify vertex $t_j$ with vertex
  $h_j^{j_3}$ and vertex $t'_j$ with vertex $f_j^{j_3}$.
\end{itemize}

\subsection{Proof of Theorem \ref{NP-complete}}
\label{proof-NP-complete}

Since the problem 3-\textsc{sat} is NP-complete, it is sufficient to
prove that the Boolean formula $F$, as defined in the previous
subsection, is true if and only if the graph $G_F$ is a facility
location graph.

Assume that the graph $G_F$ is a facility location graph and let $D$
be a directed graph such that $I(D)=G_F$.  Define an assignment of the
Boolean variables $x_i$, $i=1\dots,n$ as follows:
$$
x_i=\left\lbrace
\begin{array}{ll}
1 & \mbox{ if  the arc }g^i_1 \mbox{ enters the arc } b^i_1 \mbox{ in } D,\\
0 & \mbox{otherwise}
\end{array}
\right.
$$

Notice that from Lemma \ref{GAD1} whenever the arc $g^i_1$ enters the
arc $b^i_1$, then $g^i_j$ enters the arc $b^i_j$ for each
$j=1,\dots,m$.  Let $C_j$ be any clause of $F$. From Lemma \ref{GAD2}
(i), we must have that $r_j$ enters $r'_j$, or $s_j$ enters $s'_j$ or
that $t_j$ enters $t'_j$ in any directed graph whose intersection
graph is $\textsc{Gad}^2_j$.  We may assume that $r_j$ enters $r'_j$.
By the definition of $G_F$ the branch $r_jr'_j$ is identified with
$g^i_jb^i_j$ when $x_i$ is present in $C_j$ and in this case $x_i=1$
and so $C_j$=1. Otherwise the branch $r_jr'_j$ is identified with
$h^i_jf^i_j$ when $\bar{x}_i$ is present in $C_j$.  So the arc $h^i_j$
enters the arc $f^i_j$ and from Lemma \ref{GAD1} we have that the arc
$b^i_j$ enters the arc $g^i_j$ and by definition we have $x_i=0$,
which implies that $C_j=1$.

Now assume that there is an assignment of the variables $x_i$,
$i=1,\dots,n$ satisfying $F$.  Let us construct a directed graph $D$
such that $G_F=I(D)$. For each graph $\textsc{Gad}^1_i$ we build a
directed graph such that each arc $g^i_j$ enters the arc $b^i_j$ when
$x_i=1$ and each arc $b^i_j$ enters $g^i_j$ when $x_i=0$. This is
possible from Lemma \ref{GAD1}.  Now given a clause
$C_j=\lambda_{j_1}\vee\lambda_{j_2}\vee\lambda_{j_3}$, from Lemma
\ref{GAD2} the graph $D$ cannot exist only when the assumption (i) of
Lemma \ref{GAD2} is not satisfied. But one can check that this may
happen only when $\lambda_{j_1}=\lambda_{j_2}=\lambda_{j_3}=0$, which
is not possible.

\section{Recognizing triangle-free facility location graphs}
\label{RTFLG}

Notice that the maximum cliques on the graph $G_F$ built by our reduction in the section have size 3.  Hence it is natural to ask if recognizing
triangle-free facility location graphs remains difficult.  In this
section we show that this recognition may be done in polynomial time.

In subsection \ref{struc} we examine the structure of general FL
graphs. In subsection \ref{app-tri-free}, we restrict ourselves to
triangle-free graphs and we give the main result of this section.

\subsection{Some structural properties}
\label{struc}

In \cite{BBpg}, Ba\"{\i}ou and Barahona gave a characterization of
preimages of cycles.

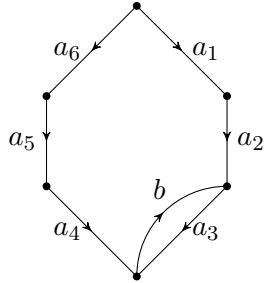
\begin{figure}[htb]
  \begin{center}
    \begin{tikzpicture}[line cap=round,line join=round,x=.6cm,y=.6cm]
      
      \coordinate (a) at (0,2);
      \coordinate (b) at ($ (a) + (2,-2) $);
      \coordinate (c) at ($ (b) + (0,-2) $);
      \coordinate (d) at ($ (c) + (-2,-2) $);
      \coordinate (e) at ($ (d) + (-2,2) $);
      \coordinate (f) at ($ (e) + (0,2) $);
      \coordinate (g) at ($(d)+(0.5,1.5)$);
      
      \fill (a) circle (1.5pt);
      \fill (b) circle (1.5pt);
      \fill (c) circle (1.5pt);
      \fill (d) circle (1.5pt);
      \fill (e) circle (1.5pt);
      \fill (f) circle (1.5pt);
      
      \begin{scope}[decoration={markings,mark=at position 0.5 with
          {\arrow{stealth'}}}]
        
        \draw[postaction={decorate}] (a) -- (b) node[midway,right]{$a_1$};
        \draw[postaction={decorate}] (b) -- (c) node[midway,right]{$a_2$};
        \draw[postaction={decorate}] (c) -- (d) node[midway,right]{$a_3$};
        
        \draw[postaction={decorate}] (e) -- (d) node[midway,left]{$a_4$};
        \draw[postaction={decorate}] (f) -- (e) node[midway,left]{$a_5$};
        \draw[postaction={decorate}] (a) -- (f)  node[midway,left]{$a_6$};
        \draw[postaction={decorate}] (d) to[bend left=45] (c);
        \draw (g) node[above]{$b$};
      \end{scope} 
    \end{tikzpicture}
  \end{center}
  \caption{The chordless cycle in $I(D)$ is
    $C=a_1,a_2,a_3,b,a_4,a_5,a_6,a_1$.  The cycle in $G$,
    $C'=a_1,\dots,a_6,a_1$.  $C''=\{b\}$.}
  \label{cycles}
\end{figure}

\begin{lemma} \cite{BBpg} 
\label{cycle}
Given a directed graph $D=(V,A)$, a subset of arcs $C\subseteq A$
induce a chordless cycle of size at least four in $I(D)$, if and only
if $C$ may be partitioned into two subsets $C'$ and $C''$, such that
$C'$ is a cycle in $D$ and there is a 1-to-1 correspondence between
the nodes in $\hat{C'}$ and the arcs in $C''$, where each arc
$(v,\bar{v})$ of $C''$ correspond to a node $v\in\hat{C'}$ where (i)
$\bar{v}\in V\setminus V(C')$ or (ii) $\bar{v}\in \tilde{C'}$ with
$\bar{v}$ is one of the two neighbors of $v$ in $C'$, see Figure
\ref{cycles} for an example.
\end{lemma}

Notice that we may have $C''=\emptyset$. In this case $C$ is a
directed cycle in $D$.  Notice the following remark.

\begin{obs}
  \label{pendent}
  Let $G$ be an undirected graph. Let $e=uv$ an edge of $G$ where $u$
  has degree one and $v$ has degree two. Then $G$ is a FL graph if
  and only if $G-u$ is also a FL graph.
\end{obs}

\begin{lemma}
  \label{sink}
  If $G$ is a FL graph, then there exists a digraph $D$ such that
  $G=I(D)$ and every sink node in $D$ has exactly one entering arc.
\end{lemma}

\proof Since $G$ is a FL graph, there exists a directed graph $D'$
such that $G=I(D')$. If $u$ is a sink node with the entering arcs
$(u_1,u),\dots,(u_k,u)$, $k\geq 2$. We may split the node $u$ and so
each arc $(u_i,u)$ is replaced by $(u_i,u_i')$, where each node $u_i'$
is a sink with only one entering arc. If $D$ is the resulting directed
graph we have that $G=I(D)$.  \qed

\begin{figure}[htb]
  \begin{center}
    \begin{tikzpicture}[line cap=round,line join=round,x=.6cm,y=.6cm]
      \coordinate (a) at (0,2);
      \coordinate (b) at ($ (a) + (0,3) $);
      \coordinate (c) at ($ (b) + (3,0) $);
      \coordinate (d) at ($ (a) + (3,0) $);
      \coordinate (e) at ($ (a) + (10,0) $);
      \coordinate (f) at ($ (e) + (0,3) $);
      \coordinate (g) at ($ (f) + (3,0) $);
      \coordinate (h) at ($ (g) + (0,-1.5) $);
      \coordinate (i) at ($ (h) + (0,-1.5) $);
      \coordinate (j) at (-2,1.5);
      \coordinate (k) at (5,1.5);
      \coordinate (l) at ($ (e) + (-2,-0.5) $);
      \coordinate (m) at ($ (i) + (+2,-0.5) $);
      
      \fill (a) circle (1.5pt);
      \fill (b) circle (1.5pt);
      \fill (c) circle (1.5pt);
      \fill (d) circle (1.5pt);
      \fill (e) circle (1.5pt);
      \fill (f) circle (1.5pt);
      \fill (g) circle (1.5pt);
      \fill (h) circle (1.5pt);
      \fill (i) circle (1.5pt);
      
      \draw [>=latex, dashed] (j) to [bend left] (k);
      \draw [>=latex, dashed] (l) to [bend left] (m);
      \draw (a) -- (b);
      \draw (b) -- (c);
      \draw (c) -- (d);
      
      \draw (e) -- (f);
      \draw (g) -- (h);
      \draw (h) -- (i);
      
      \draw (a) node[left] {$a$};
      \draw (b) node[above] {$b$};
      \draw (c) node[above] {$c$};
      \draw (d) node[right] {$d$};
      \draw (e) node[left] {$a$};
      \draw (f) node[above] {$b'$};
      \draw (g) node[above] {$b''$};
      \draw (h) node[right] {$c$};
      \draw (i) node[right] {$d$};
      
    \end{tikzpicture}
  \end{center}
  \caption{On the left: the graph $G$. On the right: the graph $G''$.}
  \label{trans}
\end{figure}
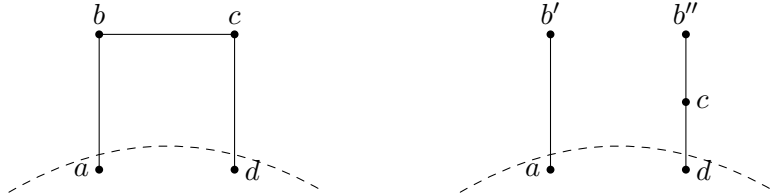

\begin{lemma} 
  \label{degre2}
  Given an undirected graph $G$. If there is an edge $e=bc$ in $G$
  such that $b$ and $c$ both have degree two and no common neighbour,
  then the following statements are equivalent:

  \begin{enumerate}
    \item[(i)] $G$ is a FL graph,
    \item[(ii)] $G - e$ is a FL graph.
  \end{enumerate}

\end{lemma}
\proof Let us call $a$ the other neighbour of $b$ and $d$ the other
neighbour of $c$ (see Figure \ref{trans}).

{\it (i) $\Rightarrow$ (ii)}. Assume that $G$ is a FL graph and hence
$G=I(D)$. Assume that the arcs $a$ and $c$ have no common vertex in
$D$, in this case the arc $b$ must share exactly one endnode with $a$
and one endnode with $c$.  Replace the arc $b=(r,s)$ by $b=(r',s)$
(respectively $b=(r,s')$) when $r$ (respectively $s$) is an endnode of $c$. The
nodes $r'$ and $s'$ are new nodes. Call the resulting graph $D'$. We
have $G - e=I(D')$.

Now assume that $a$ and $c$ have a common endnode.  Let $a=(u,v)$ and
$c=(w,t)$. Since $a$ and $c$ are not adjacent we must have $v=t$ and
$u\not=w$. Let $b=(r,s)$. If $r\not=v$, then as in the previous case
we replace the arc $b=(r,s)$ by $b=(r',s)$ (respectively $b=(r,s')$)
when $r$ (respectively $s$) is an endnode of $c$ and we obtain a new
graph $D'$ with $G - e=I(D')$. If $r=v$, replace the arc $c=(w,t)$ by
$c=(w,t')$, and if in addition we have $s=w$ we replace the arc
$b=(r,s)$ by $b=(r,s')$, where $s'$ and $t'$ are two new nodes. It is
easy to check that we obtain a new graph $D'$ with $G-e=I(D')$.

{\it (ii) $\Rightarrow$ (i).}  Now assume that $G - e=I(D')$.  Let
$G''$ be the graph obtained from $G-e$ by adding a node $b''$ and the
edge $b''c$.  To avoid confusion, we rename the node $b$ by $b'$, see
Figure \ref{trans}. From Remark \ref{pendent}, we know there exists a
directed graph $D''$ such that $G''=I(D'')$.

By Lemma \ref{sink} we may pick $D''$ such that every sink has a most
one entering arc. 

Let $b'=(r,s)$ and $b''=(t,u)$.

The nodes $b'$ and $b''$ have no node in common. Indeed, if $b'$ and
$b''$ have a common node it must be that $s=u$, and since in $G''$ the
nodes $b'$ and $b''$ are not adjacent. And $s$ must be a sink since
$b'$ and $b''$ have no common neighbour. It follows that $s$ is a sink
having at least two entering nodes, which is a contradiction.

Moreover, we may assume that $a\not=(s,r)$ and $c\not=(u,t)$. If
$a=(s,r)$. Since $a$ is the unique neighbor of $b'=(r,s)$, we may
replace the arc $b'$ by $b=(r,s')$ with $s'$ a new node. The same
arguments hold when $c=(u,t)$.

Therefore, the connections between the arcs $b'$ and $a$ and between
$b''$ and $c$ are of three types as shown in the Figure
\ref{connection} below.

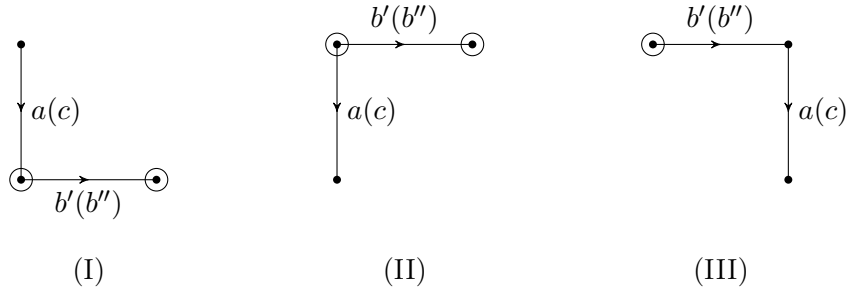
\begin{figure}[htb]
  \begin{center}
    \begin{tikzpicture}[line cap=round,line join=round,x=.6cm,y=.6cm]
      \coordinate (a) at (0,2);
      \coordinate (b) at ($ (a) + (0,-3) $);
      \coordinate (c) at ($ (b) + (3,0) $);
      \coordinate (d) at ($ (a) + (7,0) $);
      \coordinate (e) at ($ (d) + (0,-3) $);
      \coordinate (f) at ($ (d) + (3,0) $);
      \coordinate (g) at ($ (f) + (4,0) $);
      \coordinate (h) at ($ (g) + (3,0) $);
      \coordinate (i) at ($ (h) + (0,-3) $);
      \coordinate (j) at ($ (b) + (1.5,-1.5) $);
      \coordinate (k) at  ($ (e) + (1.5,-1.5) $);
      \coordinate (l) at  ($ (g) + (1.5,-4.5) $);
      
      \fill (a) circle (1.5pt);
      \fill (b) circle (1.5pt);
      \fill (c) circle (1.5pt);
      \fill (d) circle (1.5pt);
      \fill (e) circle (1.5pt);
      \fill (f) circle (1.5pt);
      \fill (g) circle (1.5pt);
      \fill (h) circle (1.5pt);
      \fill (i) circle (1.5pt);
      
      \begin{scope}[decoration={markings, mark=at position 0.5 with {\arrow{stealth'}}}]
        
        \draw (b) circle(0.25);
        \draw (c) circle(0.25);
        \draw (d) circle(0.25);
        \draw (f) circle(0.25);
        \draw (g) circle(0.25);
        \draw[postaction={decorate}] (a) -- (b) node[midway,right]{$a(c)$};
        \draw[postaction={decorate}] (b) -- (c) node[midway,below]{$b'(b'')$};
        \draw[postaction={decorate}] (d) -- (e) node[midway,right]{$a(c)$};
        \draw[postaction={decorate}] (d) -- (f) node[midway,above]{$b'(b'')$};
        \draw[postaction={decorate}] (g) -- (h) node[midway,above]{$b'(b'')$};
        \draw[postaction={decorate}] (h) -- (i)  node[midway,right]{$a(c)$};
        \draw (j) node[below] {(I)};
        \draw (k) node[below]  {(II)};
        \draw (l)  node[below]  {(III)};
      \end{scope}
    \end{tikzpicture}
  \end{center}
  \caption{The circled nodes have no other adjacency.}
  \label{connection}
\end{figure}

In some cases a simple identification of the arcs $b'$ and $b''$ give
rise to a graph $D$ with $I(D)=G$. Identifying the arcs $b'=(r,s)$ and
$b''=(t,u)$ means that we remove them and we shrink $r$ with $t$ and
$s$ with $u$, where $r'$ and $s'$ are the resulting nodes,
respectively.  Finally we put $b=(r',s')$. The cases where such an
operation may be done are summarized in the Table \ref{tab:comp}
below.

 \begin{table}[h]
  \center
  
  \begin{tabular}{cc|c|c|c|}
    \cline{3-5}
    & & \multicolumn{3}{|c|}{$c-b''$}\\
    \cline{3-5}
    & &(I)&(II)&(III)\\
    \hline
    \multicolumn{1}{|c|}{$a-b'$}& (I)& \checkmark & & \checkmark \\
    \cline{2-5}
    \multicolumn{1}{|c|}{}& (II)& & & \checkmark\\
    \cline{2-5}
    \multicolumn{1}{|c|}{}& (III)& \checkmark & \checkmark & \\
    \hline
  \end{tabular}
  \caption{Compatibility between types}
  \label{tab:comp}
\end{table}

If a connection of type (II) occurs, say between $a$ and $b'$, we may
replace the arc $a=(r,v)$ by $a=(s,v)$. The head of $a$ is unchanged
but its tail now coincide with the head of $b'$.  This transformation
will not alter $G''$, that is $G''$ is again the intersection graph of
the resulting directed graph. Moreover the connection between $a$ and
$b'$ is now of type (III). Therefore the only remaining case to study
is when the connection between $a$ and $b'$ and between $c$ and $b''$
are both of type (III). Recall that $c$ is adjacent to only $b''$ and
$d$ in $G''$.  Let $d=(v,w)$. If the head of the arc $c$ is $v$ in
$D''$, we set $c=(t,v)$.  Again $G''$ is the intersection graph of
this new directed graph and the connection of $c$ and $b''$ is of type
(II). Now assume that $w=u$, in this case we may set $b''$ to be the
arc going from the head of $c$ to $t$ and we obtain a connection of
type (I) between, $b''$ and $c$.  \qed

\medskip

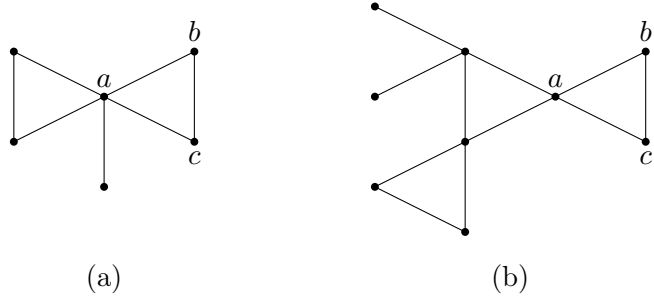
\begin{figure}[htb]
  \begin{center}
    \begin{tikzpicture}[line cap=round,line join=round,x=.6cm,y=.6cm]
      \coordinate (a) at (0,2);
      \coordinate (b) at ($ (a) + (0,-2) $);
      \coordinate (c) at ($ (b) + (2,1) $);
      \coordinate (d) at ($ (c) + (0,-2) $);
      \coordinate (e) at ($ (c) + (2,1) $);
      \coordinate (f) at ($ (e) + (0,-2) $);
      \coordinate (g) at ($ (e) + (4,1) $);
      \coordinate (h) at ($ (g) + (0,-2) $);
      \coordinate (j) at ($ (g) + (2,-1) $);
       \coordinate (k) at  ($ (j) + (0,-2) $);
       \coordinate (l) at  ($ (k) + (0,-2) $);
      \coordinate (i) at ($ (k) + (-2,-1) $);
      \coordinate (m) at ($ (j) + (2,-1) $);
      \coordinate (n) at ($ (j) + (4,0) $);
      \coordinate (o) at ($ (n) + (0,-2) $);
       \coordinate (p) at ($ (d) + (0,-1.5)$);
        \coordinate (q) at ($ (l) + (1,-0.5) $);

      \fill (a) circle (1.5pt);
      \fill (b) circle (1.5pt);
      \fill (c) circle (1.5pt);
      \fill (d) circle (1.5pt);
      \fill (e) circle (1.5pt);
      \fill (f) circle (1.5pt);
      \fill (g) circle (1.5pt);
      \fill (h) circle (1.5pt);
      \fill (i) circle (1.5pt);
      \fill (j) circle (1.5pt);
      \fill (k) circle (1.5pt);
      \fill (l) circle (1.5pt);
       \fill (m) circle (1.5pt);
      \fill (n) circle (1.5pt);
      \fill (o) circle (1.5pt);
           
        \draw (a) -- (b);
        \draw (a) -- (c);
        \draw (b) -- (c);
        \draw (c) -- (d);
        \draw (c) -- (e);
        \draw (e) -- (f);
        \draw (c) -- (f);
        \draw (j) -- (g);
        \draw (j) -- (h);
        \draw (j) -- (k);
        \draw (k) -- (l);
         \draw (k) -- (i);
         \draw (i) -- (l);
        \draw (k) -- (m);
        \draw (j) -- (m);
         \draw (m) -- (n);
         \draw (m) -- (o);
         \draw (n) -- (o);
         \draw (c) node[above]{$a$};
         \draw (e) node[above]{$b$};
         \draw (f) node[below]{$c$};
         \draw (m) node[above]{$a$};
         \draw (n) node[above]{$b$};
          \draw (o) node[below]{$c$};
          \draw (p) node[below]{(a)};
          \draw (q) node[below]{(b)};

    \end{tikzpicture}
  \end{center}
  \caption{(a):  is not a FL graph but without $e=ab$ is a FL. (b): is a FL graph, without $e=ab$ is not.}
  \label{2tri}
\end{figure}

Notice that Lemma \ref{degre2} is not true when $b$ and $c$ are two
adjacent nodes of degree two with a common neighbor $a$, that is
$a,b,c$ is a triangle. In fact, the graph shown in Figure \ref{2tri}
(a) is not a FL graph, but the graph obtained by removing the edge
$e=bc$ is a FL graph. And the graph shown in Figure \ref{2tri}(b) is a
FL graph, but if we remove the edge $e=bc$ the resulting graph is not
anymore a FL graph.  But we may easily check that when $a$ has only
one neighbor different from $b$ and $c$, then Lemma \ref{degre2}
holds.

\subsection{Application to triangle-free facility location graphs}
 \label {app-tri-free}

\begin{lemma} 
\label{fork}
Let $G$ be a connected triangle-free graph that does not contain two
adjacent nodes having both degree two. If $G$ is a FL graph, that is
there exists a directed graph $D=(V,A)$ with $G=I(D)$, then for each
vertex $v\in V$ we have (i) $v$ is the tail of at most one arc or (ii)
$v$ is the tail of exactly two arcs $(v,v_1)$ and $(v,v_2)$ and there
is no arc leaving $v_1$ and $v_2$.
\end{lemma}

\proof Since $G$ is triangle-free, $v$ is the tail of at most two
arcs. Therefore we may assume it is the tail of exactly two arcs
$(v,v_1)$ and $(v,v_2)$.  Notice that there is no arc entering $v$,
otherwise $G$ contains a triangle.  Assume that for each node $v_i$,
the arc $(v_i,v'_i)$ exists. Notice that there are no other arcs
leaving $v_1$ or $v_2$. Call $b=(v,v_1)$ and $c=(v,v_2)$. We have that
$b$ and $c$ are two adjacent nodes in $G$ having degree two.  \qed

\begin{theo}
\label{triangle-free}
Let $G$ be a triangle-free graph.  Define $G'$ to be the graph
obtained from $G$ by removing each edge $e=bc$ where $b$ and $c$ have
degree two. Then $G$ is a FL graph if and only if $G'$ admits at most
one cycle.
\end{theo}

\proof From Lemma \ref{degre2}, we may assume that $G$ does not
contain an edge $e=bc$ where $b$ and $c$ have degree two.

{\it Necessity.}  Let $G$ be a triangle-free facility location graph,
that is there exist a directed graph $D=(V,A)$, with $G=I(D)$. Also
assume that there are no two adjacent nodes of degree two. Suppose
that there is a connected component of $G$ containing two cycles $C_1$
and $C_2$.  We may assume that both $C_1$ and $C_2$ are chordless
cycles. By Lemma \ref{cycle}, $C_1$ may be partitioned into $C'_1$ and
$C''_1$ and $C_2$ may be partitioned into $C'_2$ and $C''_2$.  We have
that $C'_1$ and $C'_2$ are two cycles in $D$ and $C''_1$ and $C''_2$
are as defined in Lemma \ref{cycle}.

From Lemma \ref{fork}, we have $\dot{C'_1}=\dot{C'_2}=\emptyset$. That
is both $C'_1$ and $C'_2$ are directed cycles and that
$C''_1=C''_2=\emptyset$.  Notice that $C_1$ and $C_2$ have no node in
common, otherwise $C'_1$ and $C'_2$ must share at least one arc, and
since they are both directed cycles, we must create a triangle in
$G$. Since $C_1$ and $C_2$ belong to the same connected component we
must have a path $P$ in $G$ connecting a vertex $u$ of $C_1$ and a
vertex $v$ of $C_2$.  Let $P=(v_1=(u_1,u'_1),
v_2=(u_2,u'_2),\dots,v_{p-1}=(u_{p-1},u'_{p-1}), v_p=(u_p,u'_p))$,
with $u=v_1$ and $v=v_p$.  Then $u'_1=u'_2$ and $u'_{p-1}=u'_p$.  In
this case, $P$ cannot be a directed path in $G$. It follows, that
there exists at least one node $u_k$ that contradicts Lemma
\ref{fork}.

{\it Sufficiency.}  Consider a connected component of $G$. Suppose
that it consists of a tree. Let us construct a directed graph $D$ with
$G=I(D)$. Pick any node $r$ as a root.  Let $r=(u_0,v_0)$. Let
$r_1,\dots,r_k$ be the children of $r$ in $G$, we set $r_i=(v_i,u_0)$
for $i=1,\dots,k$.  Now each node $r_i$ play the role of $r$ and we
repeat this step.  This procedure ends with a directed graph $D$ such
that $G=I(D)$.

Suppose that there is a cycle $C$. This cycle must be chordless. Let
$C'$ be a directed cycle where each arc in $C'$ correspond to a node
in $C$. The rest of this component consist of disjoint trees each
intersect $C$ in one node. If this node is chosen to be the root of
the tree, then the procedure above may be applied to get a directed
graph $D$ such that $G=I(D)$.  \qed

As a consequence we obtain the following result.
\begin{theo} Given an undirected triangle-free graph $G=(V,E)$, we may
  decide whether or not $G$ is a facility location graph in $O(|E|)$.
\end{theo}
\proof In $O(|E|)$ we may remove all the edges $e=bc$ with both $b$
and $c$ of degree two. Then we apply a breadth-first search in
$O(|E|)$. If a node is encountered more than twice or there are two
nodes that were encountered twice, then there are two
cycles. Otherwise $G$ is a facility location graph.  \qed

\section{Consequences and related problems} 
\label{RP}

\subsection{ The vertex coloring problem}
 
A {\it vertex color} of a graph is an assignment of colors to the
nodes of the graph such that no two adjacent nodes receive the same
color.  The minimum number needed for a such coloring is called the
{\it chromatic number} and denoted by $\chi(G)$.  It is well known
that finding $\chi(G)$ is \textsc{np}-complete for triangle-free graphs.  A
direct consequence of the previous section shows that $\chi(G)\leq 3$
when $G$ is a triangle-free facility location graph.

Let $G=(V,E)$ be a triangle-free facility location problem and $G'$
the graph defined in Theorem \ref{triangle-free}.  It follows that
each connected component of $G'$ contains at most one cycle. If there
is an odd cycle then $\chi(G')=3$, otherwise $\chi(G')=2$. Let us
extend the coloring of $G'$ to $G$.  The reconstructing of $G$ from
$G'$ implies that at each step we add an edge $e=bc$ between two
pendent nodes of $G'$.  We keep calling the graph obtained at each
step $G'$, until the last step that produce $G$.  Let $b'$ and $c'$
be, respectively, the unique neighbors of $b$ and $c$ in $G'$.  If $b$
and $c$ have not the same color then we add $e=bc$ without altering
the existing coloring.  Assume that $b$ and $c$ have the same color.
If $b'$ and $c'$ have different colors, then we may assign the color
of $c'$ to $b$, or the color of $b'$ to $c$. Finally if $b'$ and $c'$
have the same color, then we pick arbitrarily $b$ or $c$ and we assign
him a third available color.

From the discussion above, we have the following:

\begin{theo} 
\label{coloration}
If $G$ is a triangle-free facility location graph, then $\chi(G)\leq
3$. Moreover, $\chi(G)$ may be computed in $O(|E|)$.
\end{theo}
A natural question arises: whether or not coloring facility locations
graphs is polynomial. Unfortunately the answer is no.  This will be
shown using a reduction from the edge coloring problem. Given an
undirected graph $G=(V,E)$ an {\it edge color} of $G$ is a coloring
of the edges that gives two different colors for each pair of incident
edges.
\begin{theo}
Coloring facility locations graphs is \textsc{np}-complete.
\end{theo}

\proof Given a undirected graph $G=(V,E)$ and positive integer $k$,
the problem of deciding wether or not we can color the edges of $G$
with $k$ colors has been proved to be \textsc{np}-complete by Holyer
in \cite{Hoy}. We will reduce this problem to the problem of whether
or not we can color the nodes of a facility location graph with $k$
colors. Precisely, from $G=(V,E)$ we build a directed graph $D$ such
that the edges of $G$ may be colored with $k$ colors if and only if
the nodes of $I(D)$ may be colored with $k$ colors.

Let $E=\{e_1,\dots,e_m\}$. We build $D$ as follows: for each edge
$e_i=uv\in E$ we add a node $v_i$ and the arcs $(u,v_i)$ and
$(v,v_i)$.  For each node $v_i$ we add $k-1$ arcs
$(v_i,v_{i_1}),\dots,(v_i,v_{i_{k-1}})$, where
$v_{i_1},\dots,v_{i_{k-1}}$ are new nodes.

Assume that $G$ admits an edge coloring with $k$ colors. Assign the
color of each edge $e_i=uv$ to the nodes $(u,v_i)$ and $(v,v_i)$ of
$I(D)$ and color the nodes $(v_i,v_{i_1}),\dots,(v_i,v_{i_{k-1}})$
with the other $k-1$ colors. This is a vertex coloring of $I(D)$ with
$k$ colors. Now assume $I(D)$ admits a vertex coloring with $k$
colors. Since each node $(v_i,v_{i_1}),\dots,(v_i,v_{i_{k-1}})$ of
$I(D)$ must receive a different colors (because they form a clique in
$I(D)$) we have that both nodes $(u,v_i)$ and $(v,v_i)$ must have the
same color, then assign this color ro the edge $e_i=uv$ of $G$. The
resulting coloring is an edge coloring of $G$ with $k$ colors or less.
\endproof
\subsection{The stable set problem}
Given an undirected graph $G=(V,E)$, a subset of nodes $S\subseteq V$
of an undirected graph is called a {\it stable set} if there is no
edge between any two nodes of $S$. The {\it maximum stable set
  problem} is to find a stable set of maximum size. This size is
usually called the {\it stability number} and denoted by $\alpha(G)$.
If we associate a weight $w(v)$ to each vertex $v\in V$, then the {\it
  maximum weighted stable set problem} if to find a stable set $S$
with $\sum_{v\in S}w(v)$ maximum.

The maximum stable set problem is \textsc{np}-complete for
triangle-free graph. One may show this result using the following
transformation due to Poljak \cite{Poljak}.  Given any undirected
graph $G=(V,E)$ replace any edge $e=uv$ in $E$ by a path
$uu',u'u'',u''v$. The resulting graph $\textsc{Sub}_G$ is
triangle-free and $\alpha(\textsc{Sub}_G)=\alpha(G)+|E|$. This shows
that the maximum stable set problem is \textsc{np}-complete in
triangle-free graphs. Using Theorem \ref{triangle-free} we have that
$\textsc{Sub}_G$ is also a facility location graph, since the removal of the
edges $u'u''$ yields a graph where each connected component is a
star. As a consequence we obtain the following result,
\begin{theo}  
\label{tfg-complete}
The maximum stable set problem is \textsc{np}-complete in
triangle-free facility location graphs.
\end{theo}
Since from Theorem \ref{coloration} one may color the vertices of any
triangle-free facility location graph with 3 colors in $O(|E|)$, this
immediately implies a 3-approximation algorithm for the maximum stable
set problem. This remains true for the maximum weighted stable set
problem. In fact, let $V_1$, $V_2$, $V_3$ be a partition of $V$ where
each subset $V_i$ is stable. Let $V'_i\subseteq V_i$, be the nodes of
$V_i$ having only positive weights, for $i=1,\dots,3$. Let
$w(V'_1)=\mbox{max }\{w(V'_2), w(V'_3)\}$ and $S^*$ the stable set of
maximum weight. We have
$$w(S^*)\leq w(V'_1)+w(V'_2)+w(V'_3)\leq 3w(V'_1).$$
\subsection{The facility location problem} 
Recall that the uncapacitated facility location problem (UFLP)
associated with a directed graph $D$ is equivalent to the maximum
weighted stable set problem with respect to $I(D)$.  Therefore, from
Theorem \ref{tfg-complete} we have the following corollary.
\begin{coro}
  The uncapacitated facility location problem associated with directed
  graph $D$ is \textsc{np}-complete even when $D$ does not contain the
  four graphs of Figure \ref{fig:triangle} as subgraphs.
\end{coro}
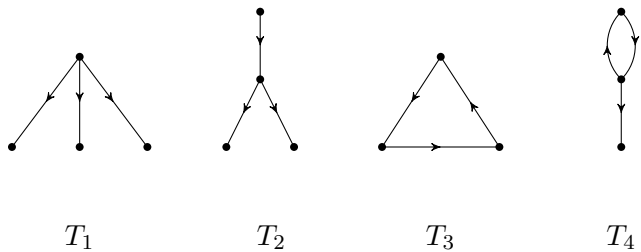
\begin{figure}[htb]
  \begin{center}
    \begin{tikzpicture}[line cap=round,line join=round,x=.6cm,y=.6cm]
      \coordinate (a) at (0,2);
      \coordinate (b) at ($ (a) + (-1.5,-2) $);
      \coordinate (c) at ($ (a) + (0,-2) $);
      \coordinate (d) at ($ (a) + (1.5,-2) $);
      \coordinate (e) at (4,1.5);
      \coordinate (f) at ($ (e) + (0,1.5) $);
      \coordinate (g) at ($ (e) + (-0.75,-1.5) $);
      \coordinate (h) at ($ (e) + (0.75,-1.5) $);
      \coordinate (k) at (8,2);
      \coordinate (i) at ($ (k) + (-1.3,-2) $);
      \coordinate (j) at ($ (k) + (1.3,-2) $);
      \coordinate (m) at (12,2);
      \coordinate (l) at ($ (m) + (0,2) $);
      \coordinate (p1) at ($ (c) + (0,-1.5) $);
       \coordinate (s1) at ($ (h) + (-0.5,-1.5) $);
       \coordinate (v1) at (8,-1.5);
       \coordinate (o) at (12,1.5);
       \coordinate (p) at ($ (o) + (0,1.5) $);
       \coordinate (q) at ($ (o) + (0,-1.5) $);
      \coordinate (o1) at  ($ (v1) + (4,0) $);
      
      \fill (a) circle (1.5pt);
      \fill (b) circle (1.5pt);
      \fill (c) circle (1.5pt);
      \fill (d) circle (1.5pt);
      \fill (e) circle (1.5pt);
      \fill (f) circle (1.5pt);
      \fill (g) circle (1.5pt);
      \fill (h) circle (1.5pt);
      \fill (i) circle (1.5pt);
      \fill (j) circle (1.5pt);
      \fill (k) circle (1.5pt);
       \fill (o) circle (1.5pt);
        \fill (p) circle (1.5pt);
         \fill (q) circle (1.5pt);

      \begin{scope}[decoration={
            markings,
            mark=at position 0.5 with {\arrow{stealth'}}}
          ] 

      \draw[postaction={decorate}] (a) -- (b);
      \draw[postaction={decorate}] (a) -- (c);
      \draw[postaction={decorate}] (a) -- (d);

      \draw (p1) node[below] {$T_1$};
      
      \draw[postaction={decorate}] (f) -- (e);
      \draw[postaction={decorate}] (e) -- (h);
      \draw[postaction={decorate}] (e) -- (g);

      \draw (s1) node[below] {$T_2$};

      \draw[postaction={decorate}] (i) -- (j);
      \draw[postaction={decorate}] (j) -- (k);
      \draw[postaction={decorate}] (k) -- (i);

       \draw (v1) node[below] {$T_3$};

      \draw[postaction={decorate}] (o) to[bend left=45] (p);
       \draw[postaction={decorate}] (p) to[bend left=45] (o);
        \draw[postaction={decorate}] (o) -- (q);
        
         \draw (o1) node[below] {$T_4$};
      
    \end{scope}
      
    \end{tikzpicture}
  \end{center}
  \caption{The forbidden subgraphs  $T_1$, $T_2$, $T_3$ and $T_4$.}
  \label{fig:triangle}
\end{figure}

In the following we will show that the UFLP remains
\textsc{np}-complete even for a more restricted class of graphs.

An undirected graph $G=(V,E)$ is called {\it cubic} is the degree of
each vertex is 3. A {\it bridge} is an edge such that its deletion
increase the number of connected components.  A {\it bridgeless} graph
is a graph with no bridge. We have the following well know result.
\begin{theo} \cite{GJS}
\label{GJ}
The maximum stable set problem in cubic graphs is \textsc{np}-complete.
\end{theo}
In \cite{GJS} it has been shown that the {\it minimum vertex cover}
problem is \textsc{np}-complete, here we look for a subset of nodes with
minimum cardinality, such that each edge has at least one endnode in
this set. Notice that if $S$ is a minimum vertex cover, then
$\bar{S}=V\setminus S$ is a maximum stable set. Then both problems
minimum vertex cover and maximum stable set are equivalent in the same
graph without any transformation.  We also notice that the proof in
\cite{GJS} use a reduction of 3-\textsc{sat} to the minimum vertex
cover problem. The graph constructed from a 3-\textsc{sat} instance is
bridgeless and each node has at most degree 3.  Moreover, each node
with degree 2 has two non-adjacent nodes of degree 3. Thus we can
remove this nodes and connect its two neighbors. It is easy to check
that if one can solve the minimum vertex cover problem in this new
graph, then one may solve it in the original graph too. From this
discussion and Theorem \cite{GJS} we have the following corollary.
\begin{coro}
\label{coro-GJS}
The maximum stable set problem in a bridgeless cubic graph is
\textsc{np}-complete.
\end{coro}
In addition to the forbidden subgraphs $T_1$, $T_2$, $T_3$ and $T_4$
we also add the subgraphs $F_1$ and $F_2$ of Figure \ref{fig:3in}, and
the UFLP remains \textsc{np}-complete.
\begin{figure}[htb]
  \begin{center}
    \begin{tikzpicture}[line cap=round,line join=round,x=.6cm,y=.6cm]
      \coordinate (a) at (0,2);
      \coordinate (b) at ($ (a) + (0,1.5) $);
      \coordinate (c) at ($ (a) + (1,1.5) $);
      \coordinate (d) at ($ (a) + (-1,1.5) $);
       \coordinate (e) at ($ (a) + (0,-1.5) $);
       
       \coordinate (a2) at ($ (a) + (4,0) $);
       \coordinate (c2) at ($ (c) + (4,0) $);
      \coordinate (d2) at ($ (d) + (4,0) $);
       \coordinate (e2) at ($ (e) + (4,0) $);

        \coordinate (s1) at ($ (e) + (0,-0.5) $);
        \coordinate (s2) at ($ (e2) + (0,-0.5) $);

      \fill (a) circle (1.5pt);
      \fill (b) circle (1.5pt);
      \fill (c) circle (1.5pt);
      \fill (d) circle (1.5pt);
      \fill (e) circle (1.5pt);
      \fill (a2) circle (1.5pt);
      \fill (c2) circle (1.5pt);
      \fill (d2) circle (1.5pt);
      \fill (e2) circle (1.5pt);

      \begin{scope}[decoration={
            markings,
            mark=at position 0.5 with {\arrow{stealth'}}}
          ]

      \draw[postaction={decorate}] (b) -- (a);
      \draw[postaction={decorate}] (c) -- (a);
      \draw[postaction={decorate}] (d) -- (a);
      \draw[postaction={decorate}] (a) -- (e);

      \draw[postaction={decorate}] (c2) -- (a2);
      \draw[postaction={decorate}] (d2) -- (a2);
      \draw[postaction={decorate}] (a2) to[bend left=45]  (e2);
       \draw[postaction={decorate}] (e2) to[bend left=45] (a2);

      \draw (s1) node[below] {$F_1$};
       \draw (s2) node[below] {$F_2$};

    \end{scope}
      
    \end{tikzpicture}
  \end{center}
  \caption{The forbidden subgraphs  $F_1$ and $F_2$.}
  \label{fig:3in}
\end{figure}
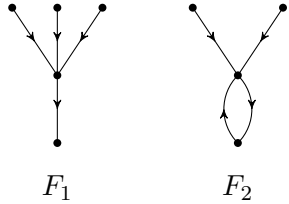
\begin{theo}
\label{UFLP-NP}
The uncapacitated facility location problem is
\textsc{np}-complete for graphs that do not contain any of $T_1$,
$T_2$, $T_3$, $T_4$, $F_1$ and $F_2$ as a subgraph.
\end{theo}
\proof Let $G=(V,E)$ be an undirected bridgeless cubic graph. From $G$
define the subdivision of it, $\textsc{Sub}_G$, as in the previous
subsection, that is each edge $e=uv\in E$ is replaced by path of size
three.  Now we construct a directed graph $D$ containing none of the
graphs $T_1$, $T_2$, $T_3$, $T_4$, $F_1$ and $F_2$ as a subgraph and
such that $I(D)=\textsc{Sub}_G$. Thus from Corollary \ref{coro-GJS}
the maximum weighted stable set problem is \textsc{np}-complete in
bridgeless cubic graphs, and by equivalence we have that UFLP is also
\textsc{np}-complete in graphs satisfying the theorem's hypothesis.
Now let us give the construction of $D$.

\begin{figure}[ht]
  \begin{center}
    \subfigure[Graph $G$]{
      \begin{tikzpicture}[line cap=round,line
        join=round,x=.6cm,y=.6cm]
      \coordinate (centre) at (0,0);
      \coordinate (ctriangle) at ($ (centre) + (-2.5,0) $);
      \coordinate (cpenta) at ($ (centre) + (2.5,0) $);
      \coordinate (t1) at ($(ctriangle) + (0:1.5)$);
      \coordinate (t2) at ($(ctriangle) + (120:1.5)$);
      \coordinate (t3) at ($(ctriangle) + (240:1.5)$);
      \coordinate (p1) at ($(cpenta) + (180:2)$);
      \coordinate (p2) at ($(cpenta) + (108:2)$);
      \coordinate (p3) at ($(cpenta) + (252:1.5)$);
      \coordinate (p4) at ($(cpenta) + (324:2)$);
      \coordinate (p5) at ($(cpenta) + (36:2)$);
      
      \draw (t1) -- (t2) -- (t3) -- (t1);
      \draw (p1) -- (p2) -- (p5) -- (p4) -- (p3) -- (p1);
      \draw (t1) -- (p1);
      \draw (t2) -- (p2);
      \draw (t3) to [bend right] (p4);
      \draw (p3) -- (p5);
      
      \fill (t1) circle (2pt);
      \fill (t2) circle (2pt);
      \fill (t3) circle (2pt);
      \fill (p1) circle (2pt);
      \fill (p2) circle (2pt);
      \fill (p3) circle (2pt);
      \fill (p4) circle (2pt);
      \fill (p5) circle (2pt);
      
    \end{tikzpicture}} \qquad  \subfigure[Graph $\textsc{Sub}_G$]{\begin{tikzpicture}[line cap=round,line join=round,x=.6cm,y=.6cm]
      \coordinate (centre) at (0,0);
      \coordinate (ctriangle) at ($ (centre) + (-2.5,0) $);
      \coordinate (cpenta) at ($ (centre) + (2.5,0) $);
      \coordinate (t1) at ($(ctriangle) + (0:1.5)$);
      \coordinate (t2) at ($(ctriangle) + (120:1.5)$);
      \coordinate (t3) at ($(ctriangle) + (240:1.5)$);
      \coordinate (p1) at ($(cpenta) + (180:2)$);
      \coordinate (p2) at ($(cpenta) + (108:2)$);
      \coordinate (p3) at ($(cpenta) + (252:1.5)$);
      \coordinate (p4) at ($(cpenta) + (324:2)$);
      \coordinate (p5) at ($(cpenta) + (36:2)$);
      
      \draw (t1) -- (t2) node[near start]{\tiny $\bullet$} node[near end]{\tiny $\bullet$};
      \draw (t2) -- (t3) node[near start]{\tiny $\bullet$} node[near end]{\tiny $\bullet$};
      \draw (t3) -- (t1) node[near start]{\tiny $\bullet$} node[near end]{\tiny $\bullet$};
      \draw (p1) -- (p2) node[near start]{\tiny $\bullet$} node[near end]{\tiny $\bullet$}; 
      \draw (p2) -- (p5) node[near start]{\tiny $\bullet$} node[near end]{\tiny $\bullet$}; 
      \draw (p5) -- (p4) node[near start]{\tiny $\bullet$} node[near end]{\tiny $\bullet$};
      \draw (p4) -- (p3) node[near start]{\tiny $\bullet$} node[near end]{\tiny $\bullet$};
      \draw (p3) -- (p1) node[near start]{\tiny $\bullet$} node[near end]{\tiny $\bullet$};
      \draw (t1) -- (p1) node[near start]{\tiny $\bullet$} node[near end]{\tiny $\bullet$};
      \draw (t2) -- (p2) node[near start]{\tiny $\bullet$} node[near end]{\tiny $\bullet$};
      \draw (t3) to [bend right] (p4) node[near start]{\tiny $\bullet$} node[near end]{\tiny $\bullet$};
      \draw (p3) -- (p5) node[near start]{\tiny $\bullet$} node[near end]{\tiny $\bullet$};
      
      \fill (t1) circle (2pt);
      \fill (t2) circle (2pt);
      \fill (t3) circle (2pt);
      \fill (p1) circle (2pt);
      \fill (p2) circle (2pt);
      \fill (p3) circle (2pt);
      \fill (p4) circle (2pt);
      \fill (p5) circle (2pt);
      
    \end{tikzpicture}}
  
\subfigure[Graph $D$]{\begin{tikzpicture}[line cap=round,line join=round,x=.6cm,y=.6cm]
      \coordinate (centre) at (0,0);
      \coordinate (ctriangle) at ($ (centre) + (-4,0) $);
      \coordinate (cpenta) at ($ (centre) + (4,0) $);
      
      \begin{scope}[decoration={
            markings,
            mark=at position 0.5 with {\arrow{stealth'}}}
          ] 

      \foreach \i in {0,1,...,8}
      {
        \fill ($(ctriangle) + ({40 * \i}:2)$) circle (1pt);
        \draw[postaction={decorate}] ($(ctriangle) + ({40 * \i}:2)$) -- ($(ctriangle) + ({40 * (\i + 1)}:2)$);
      }

      \foreach \i in {0,1,...,14}
      {
        \fill ($(cpenta) + ({180 + 24 * \i}:2)$) circle (1pt);
        \draw[postaction={decorate}] ($(cpenta) + ({180 + 24 * \i}:2)$) -- ($(cpenta) + ({180 + 24 * (\i + 1)}:2)$);
      }

      \fill (centre) circle (1pt);
      \draw[postaction={decorate}] (centre) -- ($(cpenta) + (180:2)$);
      \draw[postaction={decorate}] (centre) -- ($(ctriangle) + (0:2)$);
      
      \fill ($(centre)+(0,4)$) circle (1pt);
      \draw[postaction={decorate}] ($(centre) + (0,4)$) to [bend left] ($(cpenta) + (108:2)$);
      \draw[postaction={decorate}] ($(centre) + (0,4)$) to [bend right] ($(ctriangle) + (120:2)$);
      
      \fill ($(centre)+(0,-5)$) circle (1pt);
      \draw[postaction={decorate}] ($(centre) + (0,-5)$) to [bend right] ($(cpenta) + (324:2)$);
      \draw[postaction={decorate}] ($(centre) + (0,-5)$) to [bend left] ($(ctriangle) + (-120:2)$);

      \fill (cpenta) circle (1pt);
      \draw[postaction={decorate}] (cpenta) -- ($(cpenta) + (252:2)$);
      \draw[postaction={decorate}] (cpenta) --  ($(cpenta) + (36:2)$);
      
      \end{scope}

    \end{tikzpicture} }
  \end{center}
  \caption{From $G$ to $\textsc{Sub}_G$ to $D$}
  \label{illust}
\end{figure}
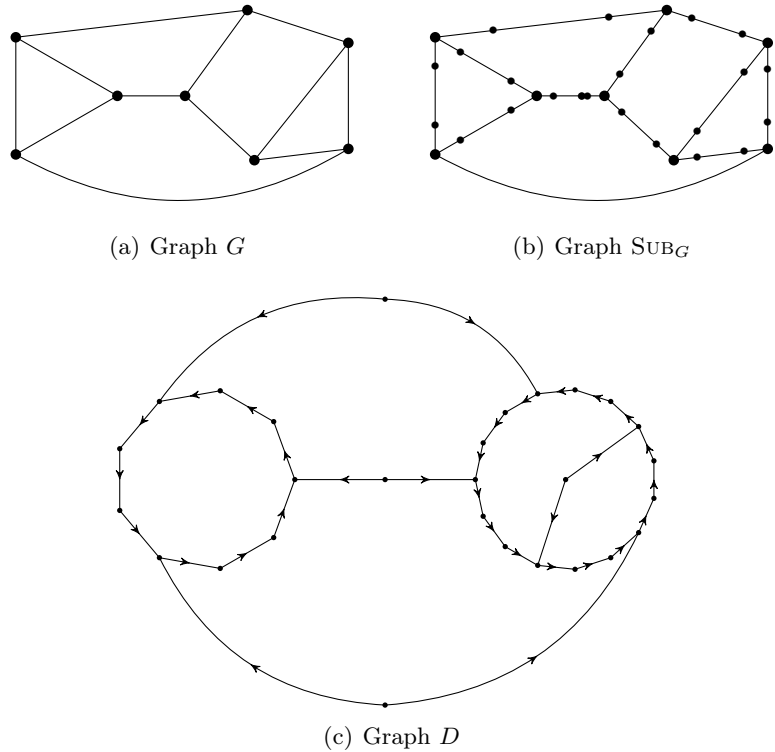

From Petersen's theorem \cite{pet}, the graph $G$ contains a perfect
matching $M$. Let $G'$ be the graph obtained by removing $M$. Each
component of $G'$ is a chordless cycle. Let $C=v_0,v_1,\dots,v_p$ be
one of these cycles. In $\textsc{Sub}_G$ this cycle corresponds to a
cycle $C'=v_0,v_1,v_2,\dots,v_{3p},v_{3p+1},v_{3p+2}$.  Let us
construct a directed graph $D$ with $I(D)=\textsc{Sub}_G$. Each cycle
$C'$ of $\textsc{Sub}_G$ may be defined in $D$ by the directed cycle
where the arc $v_i$ enters the arc $v_{i+1}$ for each
$i=0,\dots,3p+1$, and the arc $v_{3p+2}$ enters the arc $v_0$ (an arc
$a$ enters an arc $b$ means that the head of $a$ coincide with the
tail of $b$). To complete the definition of $D$ we need to consider
all the edges of $M$ and their subdivisions. Let $e=uv\in M$ and
$u_1,u_2,u_3,u_4$ the corresponding path in $\textsc{Sub}_G$.
Complete the construction of $D$ by creating for every such edge $e$
two arcs $u_2$ and $u_3$ having the same tail where $u_2$ enters the
arc $u_1$ and $u_3$ enters the arc $u_4$. This transformation is
depicted in Figure \ref{illust}.

By construction we have that $I(D)=\textsc{Sub}_G$ and that each node is the
head of at most two arcs, hence $F_1$ is not present in $D$ and it is
easy to see that with construction $F_2$ cannot occur. Also there are
no $T_1$, $T_2$, $T_3$ and $T_4$ in $D$ since $\textsc{Sub}_G$ is
triangle-free, see Figure \ref{illust}. \endproof
\section{Concluding remarks}
\label{conc}
In this paper we studied the class of facility location graphs. These
graphs come from the classical and well studied uncapacitated facility
location problem. We have shown that the recognition problem of these
graphs is \textsc{np}-complete in general and polynomially solvable in
free-triangle graphs. As a consequence, we observed that the stable
set problem still \textsc{np}-complete on a more restricted class than
triangle-free graphs and that three colors suffice to color the vertex
set of a triangle-free facility location graph.  We also studied the
complexity of two problems (1) the vertex coloring problem in facility
location graphs and (2) the uncapacitated facility location problem in
graphs that do not contain as a subgraph the graphs $T_1$, $T_2$,
$T_3$ , $T_4$, $F_1$ and $F_2$.  Let us discuss a natural attempt for
restricting more this class of graphs.

\begin{figure}[htb]
  \begin{center}
    \begin{tikzpicture}[line cap=round,line join=round,x=.6cm,y=.6cm]
      \coordinate (a) at (0,2);
       \coordinate (c) at ($ (a) + (1,1.5) $);
      \coordinate (d) at ($ (a) + (-1,1.5) $);
       \coordinate (e) at ($ (a) + (0,-1.5) $);
        
        \coordinate (a1) at ($ (a) + (4,1.5) $);
       \coordinate (c1) at ($ (c) + (4,1.5) $);
      \coordinate (d1) at ($ (d) + (4,1.5) $);
       \coordinate (e1) at ($ (e) + (4,1.5) $);
       \coordinate (f) at ($ (e1) + (0,-1.5) $);

                \coordinate (s1) at ($ (e) + (0,-0.5) $);
         \coordinate (s2) at ($ (f) + (0,-0.5) $);

      \fill (a) circle (1.5pt);
      \fill (c) circle (1.5pt);
      \fill (d) circle (1.5pt);
      \fill (e) circle (1.5pt);
      \fill (a1) circle (1.5pt);
      \fill (c1) circle (1.5pt);
      \fill (d1) circle (1.5pt);
      \fill (e1) circle (1.5pt);
      \fill (f) circle (1.5pt);

      \begin{scope}[decoration={
            markings,
            mark=at position 0.5 with {\arrow{stealth'}}}
          ]

      \draw[postaction={decorate}] (c) -- (a);
      \draw[postaction={decorate}] (d) -- (a);
      \draw[postaction={decorate}] (a) -- (e);
      
      \draw[postaction={decorate}] (c1) -- (a1);
      \draw[postaction={decorate}] (d1) -- (a1);
      \draw[postaction={decorate}] (a1) -- (e1);
      \draw[postaction={decorate}] (e1) -- (f);

       \draw (e) node[right] {$u$};
        \draw (e1) node[right] {$u$};
      \draw (s1) node[below] {$F_3$};
      \draw (s2) node[below] {$F_4$};

    \end{scope}
      
    \end{tikzpicture}
  \end{center}
  \caption{The graphs $F_3$ and $F_4$.}
  \label{Y}
\end{figure}
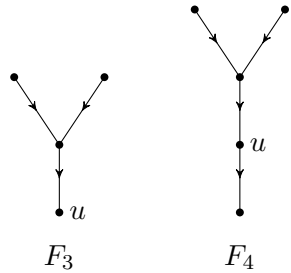

We know from \cite{BBdo1, Stauffer} that if the graph $F_3$ of Figure
\ref{Y} is forbidden, then UFLP is polynomially solvable.  Now
consider a graph without any of the subgraphs $T_1,\dots,T_4$, $F_1$
and $F_2$ and containing the subgraph $F_3$. There is at least an arc
leaving the node $u$ of $F_3$, otherwise by definition any feasible
solution of UFLP must contain $u$ and in that case $u$ may be splitted
into several copies depending on the number of arcs entering it. The
arc leaving $u$ must have a head that do not belong to $F_3$, which
lead to the graph $F_4$ of Figure \ref{Y}.  Now if we consider a
directed graph $D$ with no $F_4$ and since we do not have $F_1$ and
$F_2$, the intersection graph $I(D)$ is claw-free and hence the
maximum stable set problem is polynomially solvable \cite{Minty,
  Sbihi, Faenza}. Equivalently, the UFLP is polynomially solvable if,
in addition to the hypothesis of Theorem \ref{UFLP-NP}, we forbid also
the subgraph $F_4$.

\section*{Acknowledgements}

The authors wish to thank Reza Naserasr for fruitful discussions.

\bibliographystyle{siam}
\bibliography{paper-FLG}

\newpage

\begin{appendices}
  
  \section*{Appendix}
  
  Figures to explicitely prove Lemma \ref{GAD2}.
  
  \begin{figure}[htb]
    \begin{center}
      \subfigure[]{
        \begin{tikzpicture}[line cap=round,line join=round,x=.3cm,y=.3cm]
          \tiny
          \coordinate (x1) at (0,0);
          \coordinate (x2) at ($ (x1) + (0,-2) $);
          \coordinate (x3) at ($ (x1) + (0,-4) $);
          \coordinate (x4) at ($ (x3) + (-120:2) $);
          \node[draw,rotate=60] (x5) at ($ (x4) + (-120:2) $) {$\overline{\textsc{Inv}}$};
          \coordinate (x6) at ($ (x5) + (-120:2) $);
          \coordinate (x7) at ($ (x6) + (-120:2) $);
          \coordinate (x8) at ($ (x7) + (-150:2) $);
          \coordinate (x9) at ($ (x8) + (-150:2) $);
          \coordinate (x10) at ($ (x7) + (4,0) $);
          \coordinate (x11) at ($ (x3) + (-60:2) $);
          \node[draw,rotate=-60] (x12) at ($ (x11) + (-60:2) $) {$\overline{\textsc{Inv}}$};
          \coordinate (x13) at ($ (x12) + (-60:2) $);
          \coordinate (x14) at ($ (x13) + (-60:2) $);
          \coordinate (x15) at ($ (x14) + (-30:2) $);
          \coordinate (x16) at ($ (x15) + (-30:2) $);
          
          \begin{scope}[decoration={markings, mark=at position 0.5 with {\arrow{stealth'}}}]
            
            \draw[postaction={decorate}] (x2) -- (x1) node[midway,right]{$r'$};
            \draw[postaction={decorate}] (x3) -- (x2) node[midway,right]{$r$};
            \draw[postaction={decorate}] (x3) -- (x4) node[midway, above left]{$a$};
            \draw[postaction={decorate}] (x3) -- (x11) node[midway, above right]{$f$};
            \draw[postaction={decorate}] (x4) -- (x5.east) node[midway, above left]{$a'$};
            \draw[postaction={decorate}] (x11) -- (x12.west) node[midway, above right]{$f'$};
            \draw[postaction={decorate}] (x6) -- (x5.west) node[midway, above left]{$b'$};
            \draw[postaction={decorate}] (x13) -- (x12.east) node[midway, above right]{$e'$};
            \draw[postaction={decorate}] (x7) -- (x6) node[midway, above left]{$b$};
            \draw[postaction={decorate}] (x14) -- (x13) node[midway, above right]{$e$};
            \draw[postaction={decorate}] (x8) -- (x7) node[midway, above left]{$s$};
            \draw[postaction={decorate}] (x14) -- (x15) node[midway, above right]{$t$};
            \draw[postaction={decorate}] (x9) -- (x8) node[midway, above left]{$s'$};
            \draw[postaction={decorate}] (x15) -- (x16) node[midway, above right]{$t'$};
            \draw[postaction={decorate}] (x7) -- (x10) node[midway, below]{$c$};
            \draw[postaction={decorate}] (x10) -- (x14) node[midway, below]{$d$};
            
          \end{scope}
          
          \fill (x1) circle (1pt);
          \fill (x2) circle (1pt);
          \fill (x3) circle (1pt);
          \fill (x4) circle (1pt);
          \fill (x6) circle (1pt);
          \fill (x7) circle (1pt);
          \fill (x8) circle (1pt);
          \fill (x9) circle (1pt);
          \fill (x10) circle (1pt);
          \fill (x11) circle (1pt);
          \fill (x13) circle (1pt);
          \fill (x14) circle (1pt);
          \fill (x15) circle (1pt);
          \fill (x16) circle (1pt);

        \end{tikzpicture}
    } \quad \subfigure[]{\begin{tikzpicture}[line cap=round,line join=round,x=.3cm,y=.3cm]
        \tiny
        \coordinate (x1) at (0,0);
        \coordinate (x2) at ($ (x1) + (0,-2) $);
        \coordinate (x3) at ($ (x1) + (0,-4) $);
        \coordinate (x4) at ($ (x3) + (-120:2) $);
        \node[draw,rotate=60] (x5) at ($ (x4) + (-120:2) $) {$\overline{\textsc{Inv}}$};
        \coordinate (x6) at ($ (x5) + (-120:2) $);
        \coordinate (x7) at ($ (x6) + (-120:2) $);
        \coordinate (x8) at ($ (x7) + (-150:2) $);
        \coordinate (x9) at ($ (x8) + (-150:2) $);
        \coordinate (x10) at ($ (x7) + (4,0) $);
        \coordinate (x11) at ($ (x3) + (-60:2) $);
        \node[draw,rotate=-60] (x12) at ($ (x11) + (-60:2) $) {$\overline{\textsc{Inv}}$};
        \coordinate (x13) at ($ (x12) + (-60:2) $);
        \coordinate (x14) at ($ (x13) + (-60:2) $);
        \coordinate (x15) at ($ (x14) + (-30:2) $);
        \coordinate (x16) at ($ (x15) + (-30:2) $);
        
        \begin{scope}[decoration={markings, mark=at position 0.5 with {\arrow{stealth'}}}]
          
          \draw[postaction={decorate}] (x2) -- (x1) node[midway,right]{$r'$};
          \draw[postaction={decorate}] (x3) -- (x2) node[midway,right]{$r$};
          \draw[postaction={decorate}] (x3) -- (x4) node[midway, above left]{$a$};
          \draw[postaction={decorate}] (x3) -- (x11) node[midway, above right]{$f$};
          \draw[postaction={decorate}] (x4) -- (x5.east) node[midway, above left]{$a'$};
          \draw[postaction={decorate}] (x11) -- (x12.west) node[midway, above right]{$f'$};
          \draw[postaction={decorate}] (x6) -- (x5.west) node[midway, above left]{$b'$};
          \draw[postaction={decorate}] (x13) -- (x12.east) node[midway, above right]{$e'$};
          \draw[postaction={decorate}] (x7) -- (x6) node[midway, above left]{$b$};
          \draw[postaction={decorate}] (x14) -- (x13) node[midway, above right]{$e$};
          \draw[postaction={decorate}] (x7) -- (x8) node[midway, above left]{$s$};
          \draw[postaction={decorate}] (x14) -- (x15) node[midway, above right]{$t$};
          \draw[postaction={decorate}] (x8) -- (x9) node[midway, above left]{$s'$};
          \draw[postaction={decorate}] (x15) -- (x16) node[midway, above right]{$t'$};
          \draw[postaction={decorate}] (x7) -- (x10) node[midway, below]{$c$};
          \draw[postaction={decorate}] (x10) -- (x14) node[midway, below]{$d$};
          
        \end{scope}
        
        \fill (x1) circle (1pt);
        \fill (x2) circle (1pt);
        \fill (x3) circle (1pt);
        \fill (x4) circle (1pt);
        \fill (x6) circle (1pt);
        \fill (x7) circle (1pt);
        \fill (x8) circle (1pt);
        \fill (x9) circle (1pt);
        \fill (x10) circle (1pt);
        \fill (x11) circle (1pt);
        \fill (x13) circle (1pt);
        \fill (x14) circle (1pt);
        \fill (x15) circle (1pt);
        \fill (x16) circle (1pt);
        
      \end{tikzpicture}
    } \quad \subfigure[]{\begin{tikzpicture}[line cap=round,line join=round,x=.3cm,y=.3cm]
        \tiny
        \coordinate (x1) at (0,0);
        \coordinate (x2) at ($ (x1) + (0,-2) $);
        \coordinate (x3) at ($ (x1) + (0,-4) $);
        \coordinate (x4) at ($ (x3) + (-120:2) $);
        \node[draw,rotate=60] (x5) at ($ (x4) + (-120:2) $) {$\overline{\textsc{Inv}}$};
        \coordinate (x6) at ($ (x5) + (-120:2) $);
        \coordinate (x7) at ($ (x6) + (-120:2) $);
        \coordinate (x8) at ($ (x7) + (-150:2) $);
        \coordinate (x9) at ($ (x8) + (-150:2) $);
        \coordinate (x10) at ($ (x7) + (4,0) $);
        \coordinate (x11) at ($ (x3) + (-60:2) $);
        \node[draw,rotate=-60] (x12) at ($ (x11) + (-60:2) $) {$\overline{\textsc{Inv}}$};
        \coordinate (x13) at ($ (x12) + (-60:2) $);
        \coordinate (x14) at ($ (x13) + (-60:2) $);
        \coordinate (x15) at ($ (x14) + (-30:2) $);
        \coordinate (x16) at ($ (x15) + (-30:2) $);

        \begin{scope}[decoration={markings, mark=at position 0.5 with {\arrow{stealth'}}}]
        
        \draw[postaction={decorate}] (x1) -- (x2) node[midway,right]{$r'$};
        \draw[postaction={decorate}] (x2) -- (x3) node[midway,right]{$r$};
        \draw[postaction={decorate}] (x3) -- (x4) node[midway, above left]{$a$};
        \draw[postaction={decorate}] (x3) -- (x11) node[midway, above right]{$f$};
        \draw[postaction={decorate}] (x4) -- (x5.east) node[midway, above left]{$a'$};
        \draw[postaction={decorate}] (x11) -- (x12.west) node[midway, above right]{$f'$};
        \draw[postaction={decorate}] (x6) -- (x5.west) node[midway, above left]{$b'$};
        \draw[postaction={decorate}] (x13) -- (x12.east) node[midway, above right]{$e'$};
        \draw[postaction={decorate}] (x7) -- (x6) node[midway, above left]{$b$};
        \draw[postaction={decorate}] (x14) -- (x13) node[midway, above right]{$e$};
        \draw[postaction={decorate}] (x8) -- (x7) node[midway, above left]{$s$};
        \draw[postaction={decorate}] (x14) -- (x15) node[midway, above right]{$t$};
        \draw[postaction={decorate}] (x9) -- (x8) node[midway, above left]{$s'$};
        \draw[postaction={decorate}] (x15) -- (x16) node[midway, above right]{$t'$};
        \draw[postaction={decorate}] (x7) -- (x10) node[midway, below]{$c$};
        \draw[postaction={decorate}] (x10) -- (x14) node[midway, below]{$d$};

      \end{scope}

      \fill (x1) circle (1pt);
      \fill (x2) circle (1pt);
      \fill (x3) circle (1pt);
      \fill (x4) circle (1pt);
      \fill (x6) circle (1pt);
      \fill (x7) circle (1pt);
      \fill (x8) circle (1pt);
      \fill (x9) circle (1pt);
      \fill (x10) circle (1pt);
      \fill (x11) circle (1pt);
      \fill (x13) circle (1pt);
      \fill (x14) circle (1pt);
      \fill (x15) circle (1pt);
      \fill (x16) circle (1pt);
   
    \end{tikzpicture}
  } \quad \subfigure[]{\begin{tikzpicture}[line cap=round,line join=round,x=.3cm,y=.3cm]
      \tiny
      \coordinate (x1) at (0,0);
      \coordinate (x2) at ($ (x1) + (0,-2) $);
      \coordinate (x3) at ($ (x1) + (0,-4) $);
      \coordinate (x4) at ($ (x3) + (-120:2) $);
      \node[draw,rotate=60] (x5) at ($ (x4) + (-120:2) $) {$\overline{\textsc{Inv}}$};
      \coordinate (x6) at ($ (x5) + (-120:2) $);
      \coordinate (x7) at ($ (x6) + (-120:2) $);
      \coordinate (x8) at ($ (x7) + (-150:2) $);
      \coordinate (x9) at ($ (x8) + (-150:2) $);
      \coordinate (x10) at ($ (x7) + (4,0) $);
      \coordinate (x11) at ($ (x3) + (-60:2) $);
      \node[draw,rotate=-60] (x12) at ($ (x11) + (-60:2) $) {$\overline{\textsc{Inv}}$};
      \coordinate (x13) at ($ (x12) + (-60:2) $);
      \coordinate (x14) at ($ (x13) + (-60:2) $);
      \coordinate (x15) at ($ (x14) + (-30:2) $);
        \coordinate (x16) at ($ (x15) + (-30:2) $);

        \begin{scope}[decoration={markings, mark=at position 0.5 with {\arrow{stealth'}}}]
        
        \draw[postaction={decorate}] (x1) -- (x2) node[midway,right]{$r'$};
        \draw[postaction={decorate}] (x2) -- (x3) node[midway,right]{$r$};
        \draw[postaction={decorate}] (x3) -- (x4) node[midway, above left]{$a$};
        \draw[postaction={decorate}] (x3) -- (x11) node[midway, above right]{$f$};
        \draw[postaction={decorate}] (x4) -- (x5.east) node[midway, above left]{$a'$};
        \draw[postaction={decorate}] (x11) -- (x12.west) node[midway, above right]{$f'$};
        \draw[postaction={decorate}] (x6) -- (x5.west) node[midway, above left]{$b'$};
        \draw[postaction={decorate}] (x13) -- (x12.east) node[midway, above right]{$e'$};
        \draw[postaction={decorate}] (x7) -- (x6) node[midway, above left]{$b$};
        \draw[postaction={decorate}] (x14) -- (x13) node[midway, above right]{$e$};
        \draw[postaction={decorate}] (x7) -- (x8) node[midway, above left]{$s$};
        \draw[postaction={decorate}] (x14) -- (x15) node[midway, above right]{$t$};
        \draw[postaction={decorate}] (x8) -- (x9) node[midway, above left]{$s'$};
        \draw[postaction={decorate}] (x15) -- (x16) node[midway, above right]{$t'$};
        \draw[postaction={decorate}] (x7) -- (x10) node[midway, below]{$c$};
        \draw[postaction={decorate}] (x10) -- (x14) node[midway, below]{$d$};

      \end{scope}

      \fill (x1) circle (1pt);
      \fill (x2) circle (1pt);
      \fill (x3) circle (1pt);
      \fill (x4) circle (1pt);
      \fill (x6) circle (1pt);
      \fill (x7) circle (1pt);
      \fill (x8) circle (1pt);
      \fill (x9) circle (1pt);
      \fill (x10) circle (1pt);
      \fill (x11) circle (1pt);
      \fill (x13) circle (1pt);
      \fill (x14) circle (1pt);
      \fill (x15) circle (1pt);
      \fill (x16) circle (1pt);
      
    \end{tikzpicture}
  } \quad \subfigure[]{\begin{tikzpicture}[line cap=round,line join=round,x=.3cm,y=.3cm]
      \tiny
      \coordinate (x1) at (0,0);
      \coordinate (x2) at ($ (x1) + (0,-2) $);
      \coordinate (x3) at ($ (x1) + (0,-4) $);
      \coordinate (x4) at ($ (x3) + (-120:2) $);
      \node[draw,rotate=60] (x5) at ($ (x4) + (-120:2) $) {$\overleftrightarrow{\textsc{Inv}}$};
      \coordinate (x6) at ($ (x5) + (-120:2) $);
      \coordinate (x7) at ($ (x6) + (-120:2) $);
      \coordinate (x8) at ($ (x7) + (-150:2) $);
      \coordinate (x9) at ($ (x8) + (-150:2) $);
      \coordinate (x10) at ($ (x7) + (2,0) $);
      \coordinate (x11) at ($ (x3) + (-60:2) $);
      \node[draw,rotate=-60] (x12) at ($ (x11) + (-60:2) $) {$\overline{\textsc{Inv}}$};
      \coordinate (x13) at ($ (x12) + (-60:2) $);
      \coordinate (x14) at ($ (x13) + (-60:2) $);
      \coordinate (x15) at ($ (x14) + (-30:2) $);
        \coordinate (x16) at ($ (x15) + (-30:2) $);

        \begin{scope}[decoration={markings, mark=at position 0.5 with {\arrow{stealth'}}}]
        
        \draw[postaction={decorate}] (x2) -- (x1) node[midway,right]{$r'$};
        \draw[postaction={decorate}] (x3) -- (x2) node[midway,right]{$r$};
        \draw[postaction={decorate}] (x4) -- (x3) node[midway, above left]{$a$};
        \draw[postaction={decorate}] (x3) -- (x11) node[midway, above right]{$f$};
        \draw[postaction={decorate}] (x5.east) -- (x4) node[midway, above left]{$a'$};
        \draw[postaction={decorate}] (x11) -- (x12.west) node[midway, above right]{$f'$};
        \draw[postaction={decorate}] (x5.west) -- (x6) node[midway, above left]{$b'$};
        \draw[postaction={decorate}] (x13) -- (x12.east) node[midway, above right]{$e'$};
        \draw[postaction={decorate}] (x6) -- (x10) node[midway, above right]{$b$};
        \draw[postaction={decorate}] (x14) -- (x13) node[midway, above right]{$e$};
        \draw[postaction={decorate}] (x7) -- (x6) node[midway, above left]{$s$};
        \draw[postaction={decorate}] (x15) -- (x14) node[midway, above right]{$t$};
        \draw[postaction={decorate}] (x8) -- (x7) node[midway, above left]{$s'$};
        \draw[postaction={decorate}] (x16) -- (x15) node[midway, above right]{$t'$};
        \draw[postaction={decorate}] (x10) -- (x7) node[midway, below]{$c$};
        \draw[postaction={decorate}] (x14) -- (x10) node[midway, below]{$d$};

      \end{scope}

      \fill (x1) circle (1pt);
      \fill (x2) circle (1pt);
      \fill (x3) circle (1pt);
      \fill (x4) circle (1pt);
      \fill (x6) circle (1pt);
      \fill (x7) circle (1pt);
      \fill (x8) circle (1pt);
      \fill (x10) circle (1pt);
      \fill (x11) circle (1pt);
      \fill (x13) circle (1pt);
      \fill (x14) circle (1pt);
      \fill (x15) circle (1pt);
      \fill (x16) circle (1pt);
      
    \end{tikzpicture}
  }
\end{center}
\caption{Five possible preimages for $\textsc{Gad}^2$}
\label{fig:preimages_gad2}
\end{figure}
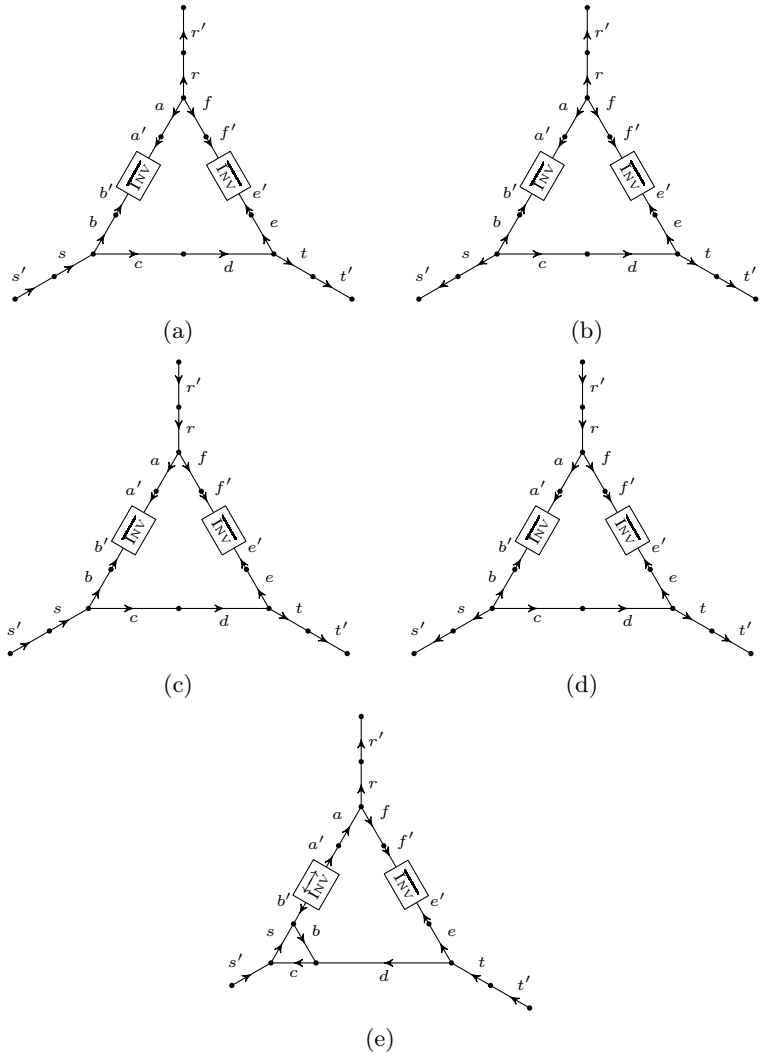

\end{appendices}

\end{document}